\documentclass[10pt]{ijnam}
\def\currentvolume{1}
\def\currentissue{1}
\def\currentyear{2004}
\def\month{March}
\pagespan{1}{18}
\copyrightinfo{2004}{} 

\usepackage{amsmath,amsthm,amssymb}
\usepackage{mathdots}
\usepackage{mathtools}
\usepackage{multicol}
\usepackage{multirow}
\usepackage{bbding}
\usepackage{subfigure}
\usepackage{xparse}
\usepackage{lipsum}
\usepackage{extarrows}
\usepackage{array}
\usepackage{courier}
\usepackage{animate}
\usepackage{dcolumn}
\usepackage[ruled,vlined,linesnumbered]{algorithm2e}
\usepackage{algorithmic}
\usepackage{pgf}
\usepackage[framemethod=TikZ]{mdframed}
\usepackage{tikz}
\usetikzlibrary{calc,shapes.multipart,chains}
\usetikzlibrary{arrows,snakes,backgrounds,shapes,patterns,shadows}
\usetikzlibrary{matrix,fit,positioning,decorations.pathmorphing}
\usepackage{color,xcolor,colortbl}

\usepackage{listings}
\tikzset{mystyle/.style n args={5}{
    rectangle,
    draw,
    fill=#1!50,
    minimum height=#2cm,
    minimum width=#3cm,
    text width=#4cm,
    align=#5
  }}

\tikzset{mystyle1/.style n args={5}{
    rectangle,
    fill=#1!50,
    minimum height=#2cm,
    minimum width=#3cm,
    text width=#4cm,
    align=#5
  }}

\tikzset{mystyle2/.style n args={5}{
    rectangle,
    fill=#1!30,
    minimum height=#2,
    minimum width=#3\textwidth,
    text width=#4\textwidth,
    align=#5
  }}

\tikzset{myarrow/.style n args={2}{
    ->,
    #1,
    line width=#2
  }}

\tikzset{list/.style={
rectangle split,
rectangle split parts=2,
rectangle split horizontal,
rectangle split part fill={red!30,blue!20},
rounded corners,
node distance=.5cm,
draw=black, thick,
minimum height=0.35cm,
text width=.5cm,
text centered,
}}

\tikzset{dlist/.style args={#1}{
rectangle split,
rectangle split parts=3,
rectangle split horizontal,
rectangle split part fill={blue!20,red!30,blue!20},
rounded corners,
node distance=#1mm,
draw=black, thick,
minimum height=0.35cm,
text centered,
}}

\tikzset{sentinel/.style args={#1}{
rectangle split,
rectangle split parts=3,
rectangle split horizontal,
rectangle split part fill={white,white,white},
rounded corners,
node distance=#1mm,
draw=black, thick,
minimum height=0.35cm,
text centered,
}}

\tikzset{square/.style args={#1}{
rectangle,draw,
fill=#1!30,
minimum height=.35cm,
minimum width=.35cm
}}

\usepackage{refcount}
\usepackage{multicol}
\newcounter{countitems}
\newcounter{nextitemizecount}
\newcommand{\setupcountitems}{%
        \stepcounter{nextitemizecount}%
        \setcounter{countitems}{0}%
        \preto\item{\stepcounter{countitems}}%
}
\makeatletter
\newcommand{\computecountitems}{%
        \edef\@currentlabel{\number\c@countitems}%
        \label{countitems@\number\numexpr\value{nextitemizecount}-1\relax}%
}
\newcommand{\nextitemizecount}{%
        \getrefnumber{countitems@\number\c@nextitemizecount}%
}
\newcommand{\previtemizecount}{%
        \getrefnumber{countitems@\number\numexpr\value{nextitemizecount}-1\relax}%
}
\makeatother    
                {\end{itemize}%
                \unskip\computecountitems\ifnumcomp{\previtemizecount}{>}{3}{\end{multicols}}{}}

\definecolor{blue}{rgb}{0.0,0.0,1.0}
\definecolor{red}{rgb}{1.0,0.0,0.0}
\definecolor{purple}{rgb}{0.75, 0.0, 1.0}
\definecolor{LightCyan}{rgb}{0.88,1,1}
\definecolor{White}{rgb}{1,1,1}

\def\cd{\cdots}

\def\vn{{\boldsymbol{n}}}

\def\vx{{\boldsymbol{x}}}

\def\vG{{\boldsymbol{G}}}

\def\vR{{\boldsymbol{R}}}
\def\vS{{\boldsymbol{S}}}
\def\vT{{\boldsymbol{T}}}

\def\vX{{\boldsymbol{X}}}

\def\vzero{{\boldsymbol{0}}}
\def\v0{{\boldsymbol{0}}}
\def\v1{{\boldsymbol{1}}}

\def\cL{{\mathcal{L}}}

\def\R{\mathbb R}

\def\vTheta{\boldsymbol{\Theta}}

\def\vxi{{\boldsymbol{\xi}}}

\def\vrho{\boldsymbol{\rho}}

\usepackage{booktabs}

\begin{document}

\title[Deep Surrogate Model for Learning Green's Function]{Deep surrogate model for learning Green's function associated with linear reaction-diffusion operator}

\author[J. Jia]{Junqing Jia}
\address{
  School of Mathematics and Statistics,
  Wuhan University,
  Wuhan, 430072, PR China
}
\email{whujjq@whu.edu.cn}


\author[L. Ju]{Lili Ju}
\address{
  Department of Mathematics,.
  University of South Carolina,
  Columbia, SC 29208, USA
}
\email{ju@math.sc.edu}

\author[X. Zhang]{Xiaoping Zhang}
\address{
  School of Mathematics and Statistics,
  Wuhan University,
  Wuhan, 430072, PR China
}
\email{xpzhang.math@whu.edu.cn}








\abstract{In this paper, we present a deep surrogate model for learning the Green's function associated with the reaction-diffusion operator in rectangular domain. The U-Net architecture is utilized to effectively capture the mapping from source to solution of the target partial differential equations (PDEs). To enable efficient training of the model without relying on labeled data, we propose a novel loss function that draws inspiration from traditional numerical methods used for solving PDEs. Furthermore, a hard encoding mechanism is employed to ensure that the predicted Green's function is perfectly matched with the boundary conditions. Based on the learned Green's function from the trained deep surrogate model, a fast solver is  developed to solve the corresponding PDEs with different sources and boundary conditions. Various numerical examples are also provided to demonstrate the effectiveness of the proposed model.
}


\maketitle

\section{Introduction}
With the rapid development and great success of deep learning technology in computer vision, natural language processing and other fields, it has also shown an increasing impact in the field of scientific computing, especially in the numerical solution of partial differential equations (PDEs) \cite{E2017DeepRitz, Sirignano2018DGM, Raissi-JCP2019}. The use of neural networks to solve PDEs has been investigated in several early works, e.g., \cite{Dissanayake94, Lagraris98}, recent advances in deep learning techniques have further stimulated new exploration in this direction.

Representative methods of interest are the physics-informed neural network (PINN) \cite{Raissi-JCP2019},  the deep Galerkin method (DGM) \cite{Sirignano2018} and the deep Ritz method (DRM) \cite{E2017DeepRitz}. All these methods model the mapping from space and/or time variables to the system states with fully connected neural network.Their differences mainly lie in the construction of loss functions. The loss functions of PINN and DGM are expressed as a weighted sum of PDE residuals at randomly selected interior points as well as solution errors at initial/boundary points.This idea also has been extended to solve inverse problems \cite{Raissi-JCP2019}, fractional differential equations \cite{Pang19}, stochastic differential equations and uncertainty qualification \cite{Yang20,Zhang19,Zhang20} and other applications. DRM \cite{E2017DeepRitz} designs loss function using the variational form of PDEs, requiring numerical integrations to train the network. Related works have subsequently emerged \cite{1912.01309,1912.03937,Wang2020,Chen2020}.

In theoretical research and engineering applications of various PDEs, including Poisson, Helmholtz and wave equations, the use of Green's function is significant. Having obtained the associated Green's function of the given differential operator, the Green's function method is used to precisely determine the solution of the corresponding PDE, which is explicitly expressed in an integral form, with the integral kernel based on the Green's function. Green's function is, in reality, a solution of the corresponding PDE with a point source subject to the homogeneous Dirichlet boundary condition. Such a problem can also be regarded as the solution of a parameterised PDE, where the location of the point source is the parameter.

However, the Green’s function in a general domain typically lacks an analytic form. Therefore, we must approximate the Green’s function numerically, which has led to increased attention on corresponding numerical methods in recent decades. Fortunately, the rapid development of deep learning techniques and their potent expressive capability potentially introduced a novel method for computing the Green's function. Supervised learning methods, such as those proposed by \cite{Gin21} and \cite{Boulle22}, have been suggested to learn the Green's functiuon. However, using these methods necessitates a considerable quantity of suitably labeled training data, which can be acquired by repeatedly solving PDEs through traditional numerical methods beforehand. The process of preparing training data consumes expensive computational resources. In addition, since these methods are purely data-driven, their generalization ability is usually restricted by the dataset coverage. In contrast, certain physics-driven models also have been proposed to compute Green's function, including GF-Net \cite{Teng22} and  BI-GreenNet \cite{Lin23}. GF-Net \cite{Teng22} extends the PINN structure \cite{Raissi-JCP2019} to solve partial differential equations stipulated by Green's function. Moreover, these models utilize certain special techniques, such as the smoothness of the Dirac delta function and domain decomposition approach to optimize the network training process. BI-GreenNet \cite{Lin23} introduces a novel framework for computing Green's function, which leverages the fundamental solution, boundary integral method and neural networks achieve high accuracy levels. 

All of the above methods are solely based on neural networks. In the past decades, traditional numerical methods, such as finite difference, finite element and finite volume methods, have been extensively studied for solving PDEs, particularly with point sources, to compute Green's function. A plausible approach is to  develop a model to compute Green's function by leveraging the benefits of both traditional methods and neural networks. In this context, we propose to use the U-Net architecture to develop a deep surrogate model for learning the Green's function of the linear reaction-diffusion operator on a rectangular domain, and to design a novel loss function, inspired by traditional numerical methods, which helps train the deep surrogate model efficiently. 

The remaining sections of the paper are organized as follow. In Section \ref{sec:problem_setting}, we briefly introduce the problem setting, including the reaction-diffusion equation, its Green's function as well as the Green's representation formula. Section \ref{sec:DSM} presents and discusses the deep surrogate model for learning the Green's function of the linear reaction-diffusion operator on a rectangular domain. This section includes the network architecture, data generation, loss function and training strategy. In Section \ref{sec:fast_solvers} we present a fast solver based on the proposed deep surrogate model to solve the corresponding PDEs. Extensive numerical experiments and comparison are provided in Section \ref{sec:numerical_experiments} to demonstrate the outstanding performance of the proposed method, including some ablation studies and the application of the deep surrogate model to the fast numerical solution of a target equation with different sources and boundary conditions.

\subsection{Problem setting and Green's function}\label{sec:problem_setting}
Let $\Omega \subset \R^d$ be a bounded Lipschitz domain, we consider the following linear reaction-diffusion operator:
\begin{equation}\label{eq:OP_RD}
  \cL(u)(\vx) := -\nabla\cdot(a(\vx)\nabla u(\vx)) + r(\vx) u(\vx), \quad \vx \in \Omega,
\end{equation}
where $a(\vx) > 0$ is the diffusion coefficient and $r(\vx)\ge 0$ is the reaction coefficient. The corresponding reaction-diffusion equation with the Dirichlet boundary condition can be represented as follows:
\begin{equation}\label{eq:PDE_RD}
  \left\{
  \begin{array}{rl}
    \cL(u)(\vx) = f(\vx), & \vx \in \Omega, \\
    u(\vx) = g(\vx), & \vx \in \partial \Omega,
  \end{array}
  \right.
\end{equation}
where $f(\vx)$ is the given source term and $g(\vx)$ gives the boundary value. The Green's function $G(\vx,\vxi)$ represents the impluse response of the PDE subject to homogenous Dirichlet boundary condition, that is, for any impulse source point $\vxi\in \Omega$, 
\begin{equation}\label{eq:PDE_Green}
  \left\{
  \begin{array}{rl}
    \cL(G)(\vx,\vxi) = \delta(\vx-\vxi), & \vx \in \Omega, \\
    G(\vx,\vxi) = 0, & \vx \in \Omega,
  \end{array}
  \right.
\end{equation}
where $\delta(\vx)$ denotes the Dirac delta source function satisfying
\begin{equation}\label{eq:delta}
  \delta(\vx) = \begin{cases}
    0, & \text{if} ~~ \vx \ne \vzero \\
    \infty, & \text{if} ~~ \vx = \vzero
  \end{cases}
  \quad\text{and}\quad
  \int_{\R^d}\delta(\vx)\,d\vx = 1. 
\end{equation}
If the Green's function $G(\vx,\vxi)$ is found, then the solution of \eqref{eq:PDE_RD} can be expressed by
\begin{equation}\label{eq:Green_representation}
  u(\vx) = \int_\Omega f(\vxi) G(\vx,\vxi)\,d\vxi - \int_{\partial\Omega} g(\vxi) a(\vxi) \frac{\partial G(\vx,\vxi)}{\partial \vn_{\vxi}}\,ds_{\vxi}, \quad \forall \vx\in \Omega, 
\end{equation}
where $\vn_{\vxi}$ denotes the unit outer normal vector on $\partial \Omega$.

\section{The deep surrogate model for learning Green's function}\label{sec:DSM}
It is noteworthy that Eq. \eqref{eq:PDE_Green} is actually a parameterized PDE with the parameter $\vxi$ and the homogeneous Dirichlet boundary conditions. We will propose a deep surrogate model to solve such a parameterized PDE, which equivalently learns the Green's function associated with the linear reaction-diffusion operator \eqref{eq:OP_RD}, and then uses it to construct a fast solver for solving the problem \eqref{eq:PDE_RD} based on the formula \eqref{eq:Green_representation}. In order to represent the Green's fucntion obeying \eqref{eq:PDE_Green}, appropriate convolutional neural network is adopted to model the mapping from the source $\delta(\vx-\vxi)$ to the solution $G(\vx,\vxi)$ of  \eqref{eq:PDE_Green}. In this work, we take the two-dimensional problem for illustration and assume $\Omega = [0,L_1]\times[0,L_2]$, but the proposed method can be naturally generalized to higher-dimensional rectangular domains.

\subsection{The U-Net architecture}

The U-Net is a representative example of a convolutional neural network (CNN), which was originally proposed for medical image segmentation, but was subsequently applied to a wide range of image processing tasks. In recent years, with the widespread application of deep learning in scientific computing, the U-Net has also been employed for regression tasks, particularly for the deep learning based method for numerical solution of PDEs, e.g. \cite{Zhao-EAAI2023}.  
Similar to all other convolutional neural networks (CNNs), the U-Net employs filter kernels for convolutional layers and pooling layers to extract features from input images. Nevertheless, the U-Net architecture is devised with a unique "U" shape, where the feature maps from the encoding path are concatenated with those of the decoding path using skip connections. This approach enables the model to capture both high-level and low-level features. Furthermore, the U-Net is recognized for its expansive path that encompasses deconvolution or up-sampling layers to progressively boost the spatial resolution of the output. 

\begin{figure}[h]
  \centering
  \includegraphics[width=\textwidth]{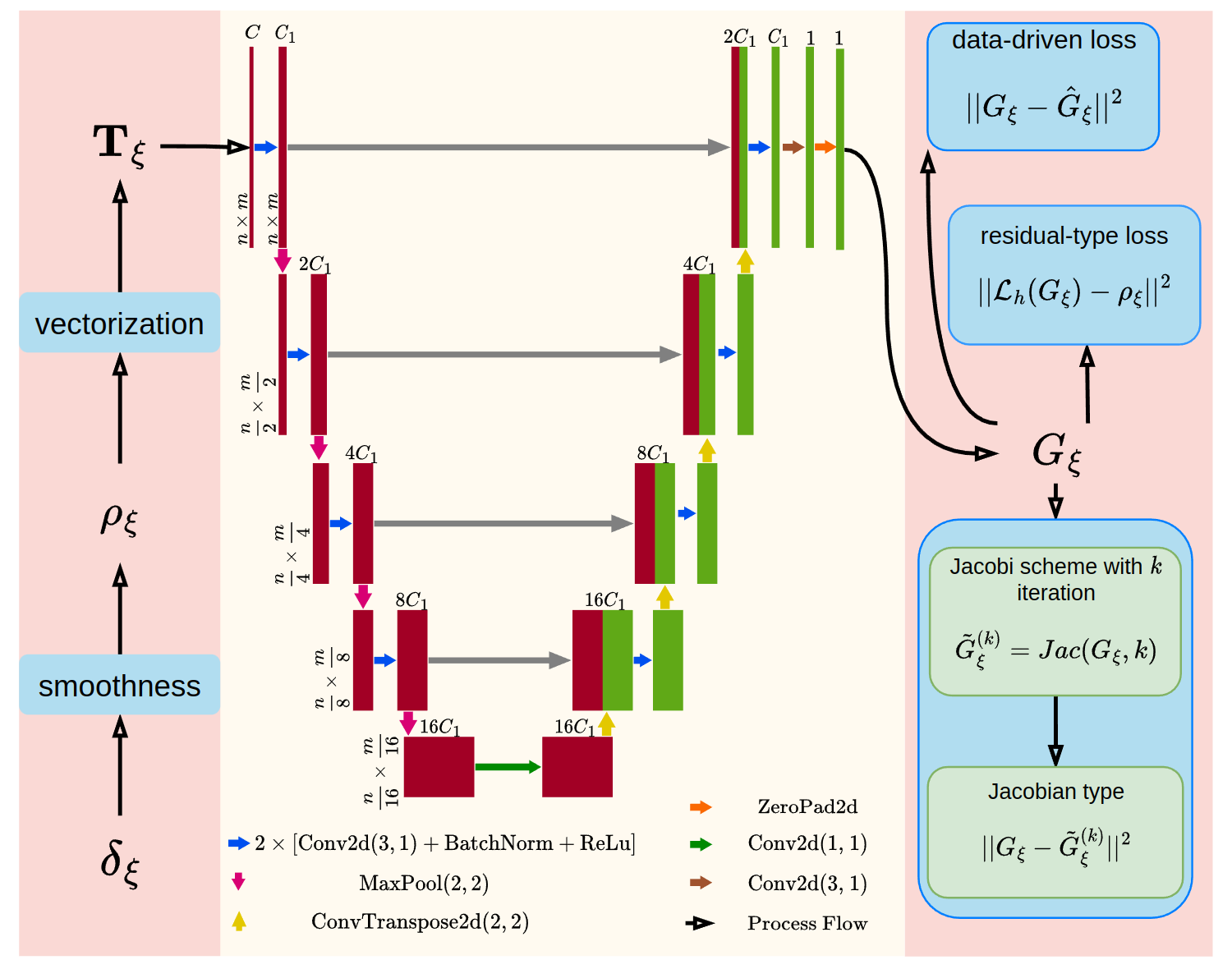}
  \caption{
    Illustration of the architecture of the U-Net for the proposed deep surrogate model for learning Green’s function. Note $\vG_{\vxi}=\mbox{\rm U-Net}(\vT_\vxi,\vTheta)$  where the input tensor $\vT_\vxi$ of the U-Net is  of dimension $n \times m \times C$ and the output tensor $\vT_\vxi$  is of dimension $n \times m$.
  }
  \label{fig:U-Net}
\end{figure}

The input tensor $\vT_\vxi$ of the U-Net is designed with dimension of $n \times m \times C$, with $C$ representing the number of input channels. The output tensor $\vG_{\vxi}$ of the U-Net is dimensioned at $n \times m$. To better suit our needs, we also slightly modify the classic architecture of the U-Net by introducing two hyper-parameters. One of these is the channels of the first hidden layer, denoted as $C_1$, which identifies the number of extracted features in the begining. The other one is the depth of the encoder/decoder, denoted as $D$. 
As depicted in Figure \ref{fig:U-Net}, each encoding operation in the U-Net downsamples the input size of the previous layer while simultaneously doubling the channel number of the input tensor. Conversely, each decoding operation in the U-Net doubles the input size of the previous layer and halves the channel number. By adding more coding and decoding layers, the depth of this architecture can be easily increased. To realize the hard encoding of the homogenous boundary condition  obeyed by Green's function, we add a zero padding operation at the end of the architecture.

\subsection{Approximation of the Delta delta function}

The Dirac delta function \eqref{eq:delta} is approximated by a multidimentional Gaussian density function
\begin{equation}\label{eq:delta_approx}
  \delta(\vx-\vxi)\approx \rho(\vx-\vxi) = \frac1{\left(\sqrt{2\pi} \sigma\right)^2} \exp\left(-\frac{|\vx-\vxi|^2}{2\sigma^2}\right),
\end{equation}
where the parameter $\sigma>0$ denotes the standard deviation of the distribution. As $\sigma\to 0$, the function \eqref{eq:delta_approx} converges to the Dirac delta function pointwisely except at the point $\vx=\vxi$. In practice, the standard deviation $\sigma$ is set to be a value proportional to the mesh size of the problem domain $\Omega$. 

\subsection{Data generation}\label{sec:data_gen} 
Let us uniformly partition the domain $\Omega$ in each direction to
obtain a rectangular mesh with nodes $\vX = \{\vx_{i,j}=((i-1)h_1,(j-1)h_2)\,|\, i=1,\cdots,n,\,j=1,\cdots,m\}$, where 
$h_1= L_1/(n-1)$ and $h_2= L_2/(m-1)$.
 For each fixed source point $\vxi$,  we first compute the distance between $\vxi$ and each node in the mesh $\vX$, and then assemble them into an array 
\begin{equation}\label{eq:R}
  \vR_{\vxi} = \left\{\|\vx_{i,j}-\vxi\|_2 ~|~ \vx_{i,j} \in \vX\right\}.
\end{equation}
Next a normalization is adopted to get  
\begin{equation}
  \vR_{\vxi} \leftarrow \frac{\vR_{\vxi} - R_{\min}}{R_{\max} - R_{\min}},
\end{equation}
where $R_{\min}$ and $R_{\max}$ are the minimum and maximum of $\vR_\vxi$, respectively. The right-hand term $\delta(\vx-\vxi)$ of \eqref{eq:PDE_Green} is also  evaluated on  $\vX$ for each  $\vxi$, which leads to an array 
\begin{equation}\label{eq:Gaussian}
  \vrho_\vxi = \left\{\rho(\vx_{i,j}-\vxi) ~|~ \vx_{i,j} \in \vX\right\},
\end{equation}
where $\rho(\cdot)$ is defined by \eqref{eq:delta_approx}.

Since the source point $\vxi$ can be randomly sampled at any location in the solution domain $\Omega$, we can easily generate the samples.  In our experiments, $2000$ training samples and $100$ validation samples are generated  with the uniform distribution for $\vxi$ in $\Omega$.
The  input tensors \{$\vT_\vxi$\} consists of  three types, including $1$-channel input $\vT_\vxi^{(1)}=\vrho_\vxi$,  $2$-channels input $\vT_\vxi^{(2)}=[\vR_\vxi, \vrho_\vxi]$, and  $3$-channels input $\vT_\vxi^{(3)} = [\vX, \vrho_\vxi]$. 


\subsection{Loss function}\label{sec:loss}
To train the deep surrogate model in a physics-driven fashion, we need to construct a loss function based on the PDE \eqref{eq:PDE_Green}.
Unlike PINN and its variations, we will not use the  strong form of the PDE. Instead, we discretize \eqref{eq:PDE_Green} by conventional numerical schemes. Specifically, we adopt the second-order central finite difference scheme to discretize \eqref{eq:PDE_Green} on $\vX$, which leads to 
\begin{equation}\label{eq:FD_Green}
  \cL_h(G_h)(\vx_{i,j}, \vxi) = \rho(\vx_{i,j}-\vxi). 
\end{equation}
where  $\vG_{\vxi}=\{G_h(\vx_{i,j}, \vxi)~|~\vx_{i,j} \in \vX\}$ and $\cL_h$ is the discrete operator for approximation of the differential operator $\cL$ given as follows:
\begin{equation}
\begin{aligned}
\cL_h(G_h)(\vx_{i,j},\vxi)  =&\;c_{i,j}G_h(\vx_{i,j}, \vxi)-c_{i+1,j}G_h(\vx_{i+1,j}, \vxi)\\
&\;-c_{i-1,j}G_h(\vx_{i-1,j}, \vxi)-c_{i,j+1}G_h(\vx_{i,j+1}, \vxi)\\
&\;-c_{i,j-1}G_h(\vx_{i,j-1}, \vxi)
+r_{ij}G_h(\vx_{i,j}, \vxi).
\end{aligned}
\end{equation}
where $r_{ij}=r(\vx_{i,j})$ and
$$c_{i+1,j} = a(\vx_{i+1/2,j})/h_1^2, \quad c_{i-1,j} = a(\vx_{i-1/2,j})/h_1^2, $$ 
$$c_{i,j+1} = a(\vx_{i,j+1/2})/h_2^2,\quad c_{i,j-1} = a(\vx_{i,j-1/2})/h_2^2,$$
$$c_{i,j} = c_{i+1,j} +c_{i-1,j} +c_{i+1,j+1} +c_{i,j-1}.$$

Then a  natural  and common way  to construct the loss function  is to use the residual of \eqref{eq:FD_Green}:
\begin{equation}\label{eq:res_loss}
\begin{aligned}
  Loss_{\text{res}}(\vTheta) =& \sum_{\vxi} \|\cL_h(\vG_{\vxi}) - \vS_\vxi \|^2 \\
  =& \sum_{\vxi} \sum_{i,j}|\cL_h(G_h)(\vx_{i,j},\vxi) - \rho(\vx_{i,j}-\vxi) |^2,
\end{aligned}
\end{equation}
where $ \vS_\vxi  = \{\rho(\vx_{i,j}-\vxi)~|~\vx_{i,j} \in \vX\}$,
which is referred as the \textbf{residual-type loss}. It is a discrete analogue of the loss function commonly used in PINN. Unfortunately, numerical experiments in Section \ref{sec:numerical_experiments} exhibit that the use of such loss function is quite hard to train the proposed deep surrogate model and could lead to a poor performance.

Inspired by the idea of Jacobi  iterative scheme for solving linear systems, we  propose and test a new loss function  defined by
\begin{equation}\label{eq:jac_loss}
  Loss_{\text{jac}}(\vTheta) = \sum_{\vxi} \|\vG_\vxi - \tilde \vG_\vxi^{(k)} \|^2,
\end{equation}
where $\tilde\vG_\vxi^{(k)}$ is the approximate solution of \eqref{eq:FD_Green} obtained by using Jacobi iteration scheme with the initial value $\tilde \vG_\vxi^{(0)} = \vG_\vxi$ and $k$ iterations, i.e., 
\begin{equation}
\begin{aligned}
\tilde G_h^{(l+1)}(\vx_{ij},\vxi) =&\;\dfrac{1}{c_{i,j}+r_{ij}}\Big[\rho(\vx_{i, j}-\vxi) +c_{i+1,j}\tilde G_h^{(l)}(\vx_{i+1,j}, \vxi))\\
&\;\qquad+c_{i-1,j}\tilde G_h^{(l)}(\vx_{i-1,j}, \vxi))+c_{i,j+1}\tilde G_h^{(l)}(\vx_{i,j+1}, \vxi))\\
&\;\qquad+c_{i,j-1}\tilde G_h^{(l)}(\vx_{i,j-1}, \vxi))
\Big], \qquad l = 0,1,\cd,k-1.
\end{aligned}
\end{equation}
We will refer \eqref{eq:jac_loss} as the  \textbf{Jacobi-type loss}.

For comparison purposes, we also consider and test a {\bf data-driven} loss function as follows:
\begin{equation}\label{eq:datadriven_loss}
  Loss_{\text{data}}(\vTheta) = \sum_{\vxi} \|\vG_\vxi - \hat \vG_\vxi \|^2,
\end{equation}
where $\hat \vG_\vxi$ is obtained by taking the final convergent result of the Jacobi  iterative solution $\tilde \vG_\vxi^{(k)}$, i.e., $\hat \vG_\vxi = \lim\limits_{k\rightarrow\infty} \tilde \vG_\vxi^{(k)}$.

\subsection{Training strategies}
This section explores training strategies for the deep surrogate model equipped with $Loss_{\text{jac}}$. The objective of the training process is to form a virtuous circle through gradually optimizing the network from the approximate solutions generated by the Jacobi iteration method. The U-Net's predictions can then be served as a potentially improved initial solutions for the Jacobi iteration in the subsequent training step.

Three options for choosing the optimal iteration number  $k$ in \eqref{eq:jac_loss} are considered. Using a fixed $k$ in the Jacobi iteration scheme during the training process is a conventional approach, referred ``constant strategy''.  In this approach selecting an optimal $k$ is important in order to  balance accuracy and computational complexity. The second approach is to first set a larger value for $k$ and then gradually decreasing it as the training progresses until it reaches a small value, which is  referred  as the ``dynamic strategy''.   A more reasonable approach is to adaptively adjust $k$ by comparing the validation errors observed in two successive epochs. If the error observed in the current epoch is significantly greater than that of the previous epoch, then $k$ should be increased, and conversely, if it is smaller then $k$ needs to be decreased. This approach is referred as ``adaptive strategy''. 


\section{Fast PDE solver based on the learned Green's function}\label{sec:fast_solvers}
Once the deep surrogate model is trained, numerical solution of the linear reaction-diffusion problem \eqref{eq:PDE_RD} can be directly computed based on the Green's formula \eqref{eq:Green_representation} through the learned Green's function. To ensure accurate evaluation of the integrals in  \eqref{eq:Green_representation} accurately, we apply numerical quadrature on rectangular meshes. To achieve this, we use the rectangular mesh of the domain $\Omega$ for training the deep surrogate model, which consists of rectangles $\mathcal R_q = \{R_l\}$. Let us denote the intersection of the rectangle edges with the domain boundary by $\mathcal E_q^{\text{bdry}}=\{E_m\}$.  By using the symmetry of Green's function, we have 
\begin{equation}\label{eq:Quad}
\begin{aligned}
  u(\vxi)
  &\approx \sum_{R_l \in \mathcal R_q} I_{\vx,h}^{R_l} [f(\vx)G(\vx,\vxi)] \\
&\qquad - \sum_{E_m \in \mathcal E_q^{\text{bdry}}} I_{\vx,h}^{E_m} [g(\vx) a(\vx) (\nabla_\vx G(\vx,\vxi) \cdot \vn_\vx)],
\end{aligned}
\end{equation}
where $I_{\vx,h}^{R_l} [\cdot]$ denotes the numerical quadrature for evaluating $$\int_{R_l} f(\vx)G(\vx,\vxi)\,d\vx$$ and $I_{\vx,h}^{E_m} [\cdot]$ the numerical quadrature for evaluating $$\int_{E_m} g(\vxi) a(\vxi) (\nabla_\vx G(\vx,\vxi) \cdot \vn_\vx)\,d s_\vx,$$ respectively.

\section{Numerical experiments} \label{sec:numerical_experiments}

This section presents various numerical experiments. We first conduct ablation studies for the deep surrogate model used to learn the Green's function of the Laplacian operator. Then, we test more examples on the learned Green's functions of the reaction-diffusion operator and corresponding fast solver.
In the following examples, the solution domain is chosen to be $[-1,1] \times [-1,1]$ and partitioned into a uniform rectangular mesh of $64\times 64$ uniform nodes, i.e., $n=m=64$. In all experiments, the maximum number of epochs and the batch size are to $150$ and $6$, respectively. All experiments were implemented using the PyTorch framework and run on the GTX 2080Ti cards.

\subsection{Ablation study of deep surrogate model} 
To simplify the matter, we use the deep surrogate model for learning the Green's function of the Laplacian operator  (i.e., $a(\vx) \equiv 1$ and $r(\vx) \equiv 0$) as an example. We  conduct a series of ablation studies to measure the influence of the model's performance, including the impact of network architecture, loss functions, input forms and the number of Jacobi  iterations. In this subsection, the number of Jacobi iterations remains fixed at $k=20$ (constant strategy) for $Loss_{\text{jac}}$ except for the experiments in subsections \ref{sec:Jacobi_fixed} and \ref{sec:Jacobi_dynamic}. 
 
\subsubsection{Effect of the U-Net architecture}\label{sec:UNet}
The U-Net architecture used for the proposed deep surrogate model is determined by the number of channels of the first hidden layer ($C_1$) and the depth of its encoder/decoder ($D$), as already explained in subsection 2.1. We carefully investigate its effect on the performance of the model, and report the corresponding test results on the model sizes and the three training MSE losses (i.e., the residual-type loss $Loss_{\text{res}}$, the Jacobi-type loss $Loss_{\text{jac}}$, and the data-driven loss $Loss_{\text{data}}$) for the U-Net architecture under various values of $C_1$ and $D$ in Table \ref{tab:1}. Our observations include: 1)  the prediction of the model equipped with the residual-type loss ($Loss_{\text{res}}$) is always unsatisfactory regardless of the choice of the U-Net architecture; 2)  for a fixed depth $D$, the performance of the model will gradually improve as the number of channels $C_1$ increases; 3)  for the model equipped with the Jacobi-type loss $Loss_{\text{jac}}$, the performance improvement of the model does not continue when the depth $D$ increases up to a certain level. To balance the size and performance of the proposed deep surrogate model, we will use the U-Net architecture with $C_1=32$  and $D=4$  in the subsequent analysis, which appears to perform the best in all cases based on Table \ref{tab:1}.

\begin{table}[htbp]
  \centering
  \begin{tabular}{p{1.5cm}<{\centering}p{1.5cm}<{\centering}p{2cm}<{\centering}p{1.5cm}<{\centering}p{1.5cm}<{\centering}p{1.5cm}<{\centering}}
    \toprule
    $C_1$ & $D$ & Model Size & $Loss_{\text{res}}$ & $Loss_{\text{jac}}$ & $Loss_{\text{data}}$\\\midrule
    4  &   3 &   15.1K &3.69e-3 &  2.29e-5 &  1.85e-4\\
    8  &   3 &   59.1K &2.60e-3 &  1.16e-5 &  7.93e-6\\
    16 &   3 &   234K  &3.10e-3 &  8.25e-6 &  3.41e-6\\
    32 &   3 &   930K  &2.88e-3 &  7.79e-6 &   2.65e-6\\\midrule
    4  &   4 &   59.7K &3.25e-3 &  2.39e-6 &  8.09e-6\\
    8  &   4 &   236K  &2.67e-3 &  1.49e-6 &  2.16e-6\\
    16 &   4 &   940K  &2.90e-3 &  1.61e-6 &  1.40e-6\\
    32 &   4 &   3.8M  &3.56e-3 &  \underline{1.15e-6} &  1.67e-6\\\midrule
    4  &   5 &   237K  &3.83e-3 &  2.58e-6 &   6.74e-6\\
    8  &   5 &   943K  &3.33e-3 &  6.85e-6 &  2.10e-6\\
    16 &   5 &   3.8M  &2.86e-3 &  1.61e-6 &  1.49e-6\\
    32 &   5 &    15M  &2.52e-3 &  1.22e-6 &  1.16e-6\\
    \bottomrule
  \end{tabular}\vspace{0.2cm}
  \caption{Results on the model sizes and the training MSE losses  for the U-Net architecture under various values of $C_1$ and $D$. }
  \label{tab:1}
\end{table}

\subsubsection{Effect of the  loss functions}\label{sec:Loss}

\begin{figure}[htbp]
  \centering
  \begin{tabular}{ccc}
    \includegraphics[width=0.48\textwidth]{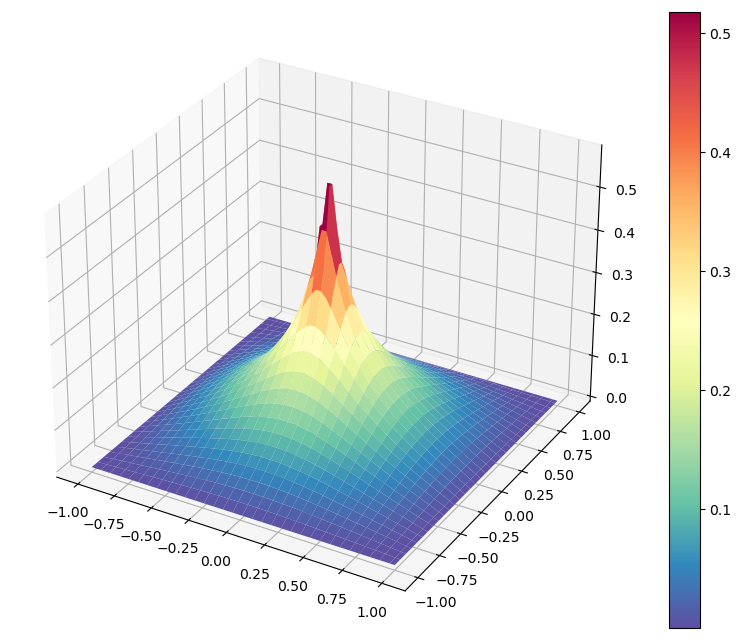} &
    \includegraphics[width=0.48\textwidth]{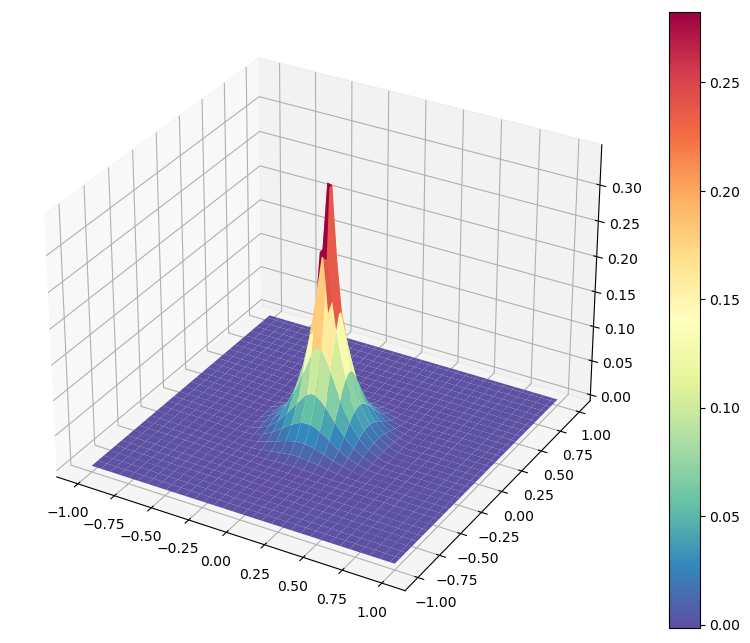} \\
    (a) Reference  & (b) $Loss_{\text{res}}$ \\
    \includegraphics[width=0.48\textwidth]{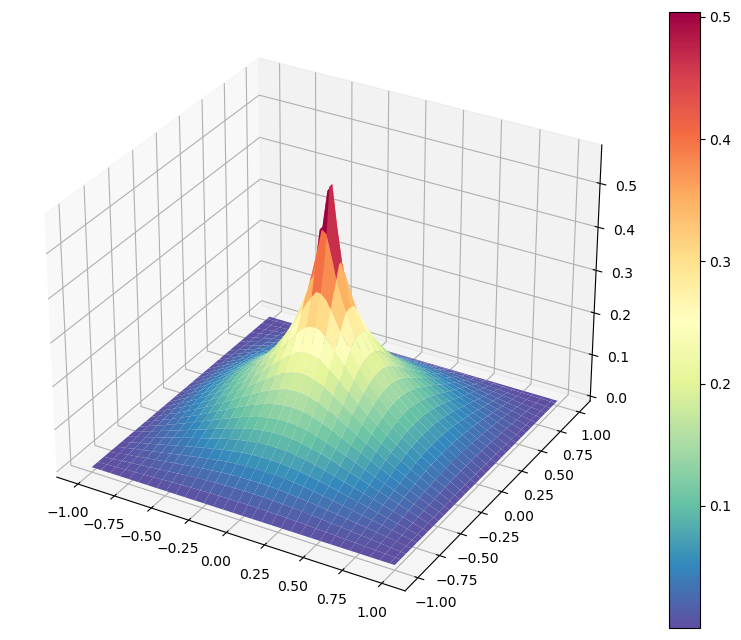} &
    \includegraphics[width=0.48\textwidth]{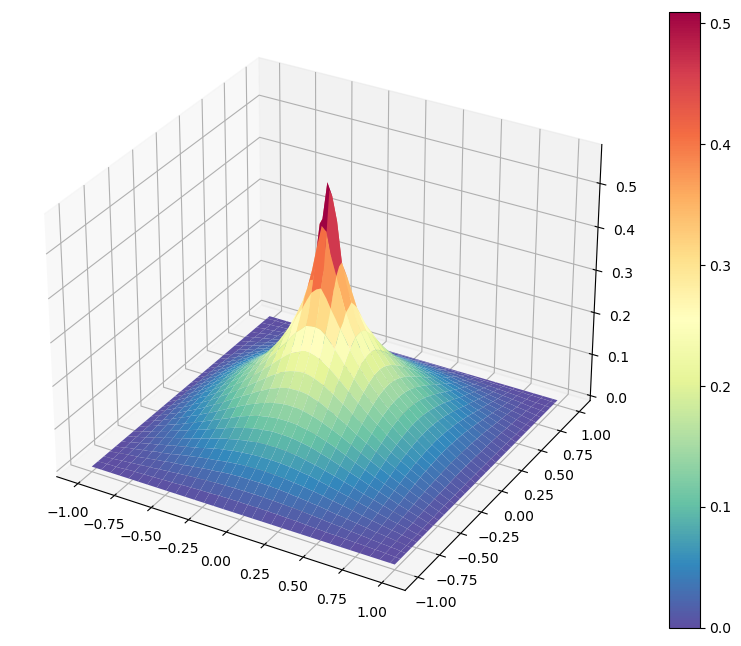} \\
    (c) $Loss_{\text{jac}}$  & (d) $Loss_{\text{data}}$ \\
  \end{tabular}
  \caption{Plots of Green's functions $G(\vx,\vxi)$ at $\vxi=(0,0)$ computed by the finite difference method (as the reference solution) and predicted by the proposed deep surrogate model equipped with the three different loss functions ($Loss_{\text{res}}$, $Loss_{\text{res}}$, $Loss_{\text{data}}$) respectively. }
  \label{fig:Surface_Loss}
\end{figure}

The key of training the proposed deep surrogate model often lies in the choice of loss functions.  Figure \ref{fig:Surface_Loss} present a visual comparison of the 
the Green's function computed by the finite difference method (as the reference solution) with  those predicted by the proposed deep surrogate model equipped with the three different loss functions.  The following observations are made: 1) the model equipped with the residual-type loss learns the rough shape of the Green's function but its detailed values are almost completely inaccurate; 2) the results predicted by the model equipped with the Jacobi-type loss and data-driven loss are very similar and both are quite accurate.

\begin{figure}[htbp]
  \centering
  \begin{tabular}{ccc}
    $e_2 = 2.23\times 10^{-3}$& $e_2 = 1.90\times 10^{-3}$\\
    \includegraphics[width=0.43\textwidth]{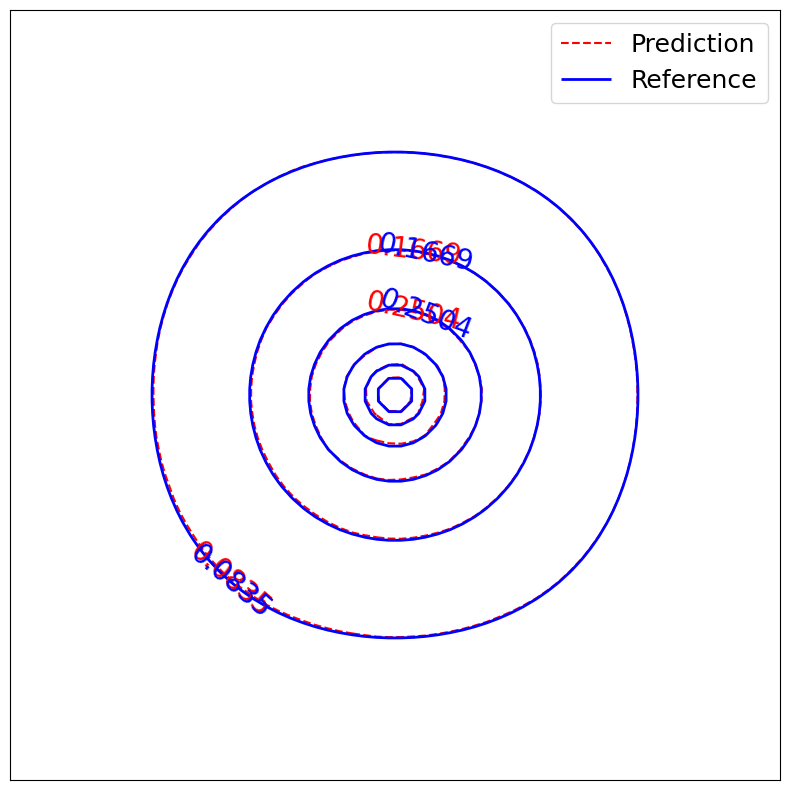} &
    \includegraphics[width=0.43\textwidth]{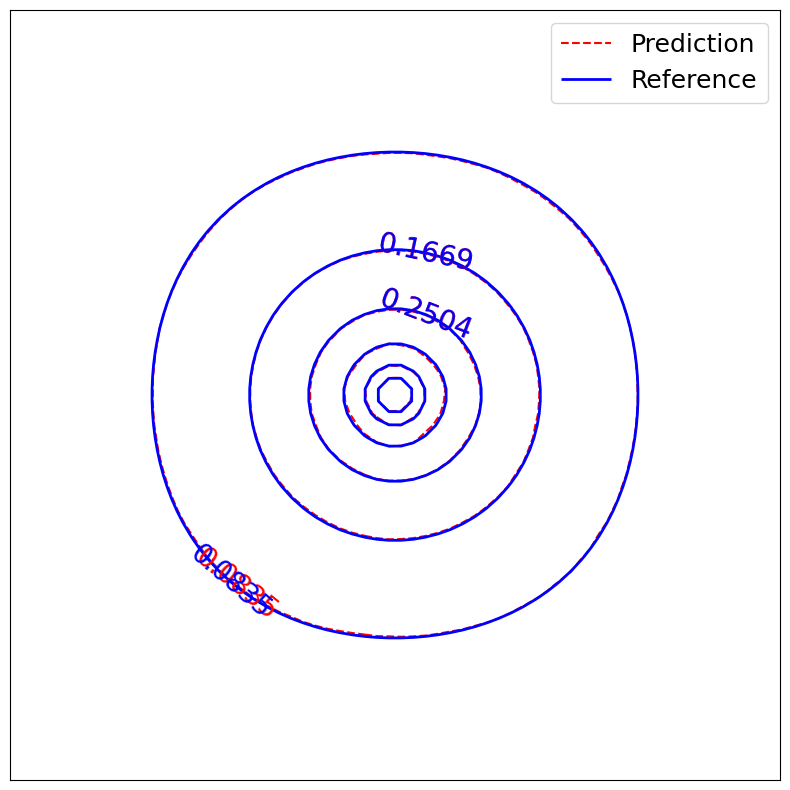} \\
    (a) $Loss_{\text{jac}}$  & (b) $Loss_{\text{data}}$ \\
  \end{tabular}
  \caption{Comparisons of the contour maps of Green's functions $G(\vx,\vxi)$ at $\vxi=(0,0)$ between the reference solution and the predicted solutions by the proposed deep surrogate model equipped with $Loss_{\text{jac}}$  and $Loss_{\text{data}}$ respectively. The corresponding $L^2$ errors (denoted as $e_2$) are also provided.}
  \label{fig:3}
\end{figure}

\begin{figure}[htbp]
  \centering
  \begin{tabular}{ccc}
    \includegraphics[width=0.43\textwidth]{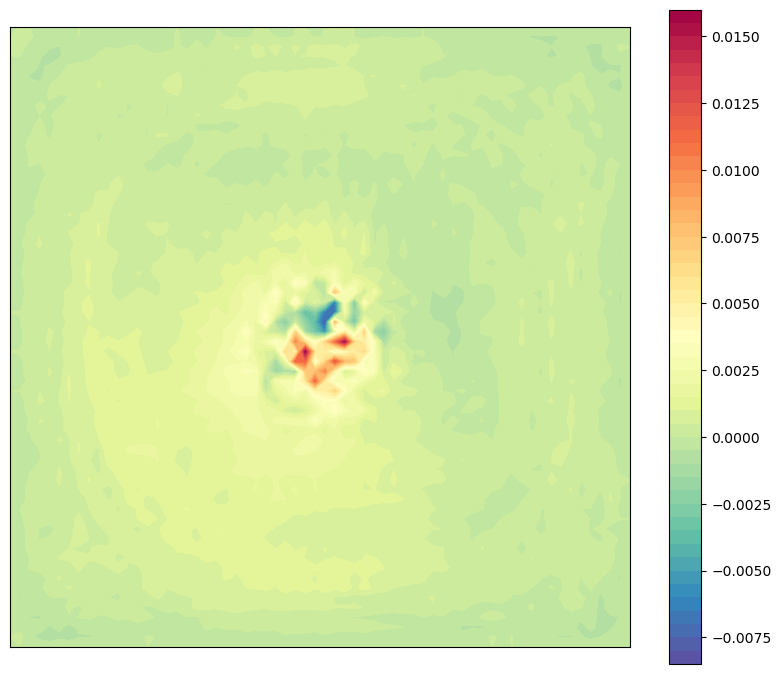} &
    \includegraphics[width=0.43\textwidth]{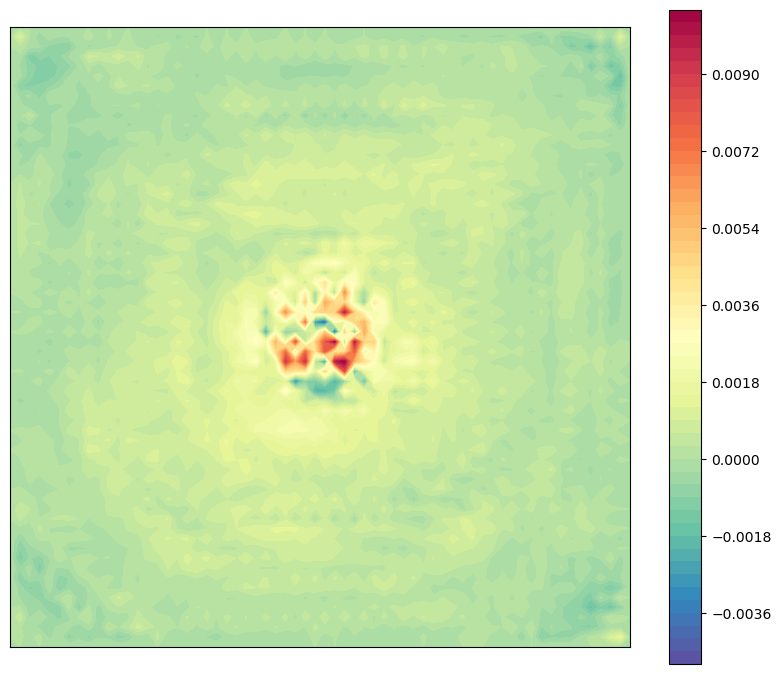} \\
    (a) $Loss_{\text{jac}}$  & (b) $Loss_{\text{data}}$ \\
  \end{tabular}
  \caption{Heat maps of the errors of Green's functions $G(\vx,\vxi)$ at $\vxi=(0,0)$ for the predicted solutions by using $Loss_{\text{jac}}$ and $Loss_{\text{data}}$ respectively.}
  \label{fig:4}
\end{figure}

Comparisons of the contour maps of the reference solution and the predicted solutions with $Loss_{\text{jac}}$ and $Loss_{\text{data}}$ are provided  in Figure \ref{fig:3},
together with the corresponding $L^2$ errors. It is observed that for the the models equipped with $Loss_{\text{jac}}$ and $Loss_{\text{data}}$, the contour lines (indicating the gradient information) of the predicted solutions overlap well with those of the reference solution. Figure \ref{fig:4} presents the heat maps of the errors for the predicted solutions  by using $Loss_{\text{jac}}$ and $Loss_{\text{data}}$ , from which, we find that the predictive errors of the proposed model equipped with $Loss_{\text{jac}}$ are comparable to that of the model equipped with  $Loss_{\text{data}}$.

\subsubsection{Effect of input forms}\label{sec:Input}

\begin{figure}[htbp]
  \centering
  \begin{tabular}{cc}
    & $e_2 = 4.33 \times 10^{-3}$\\
    \includegraphics[width=0.43\textwidth]{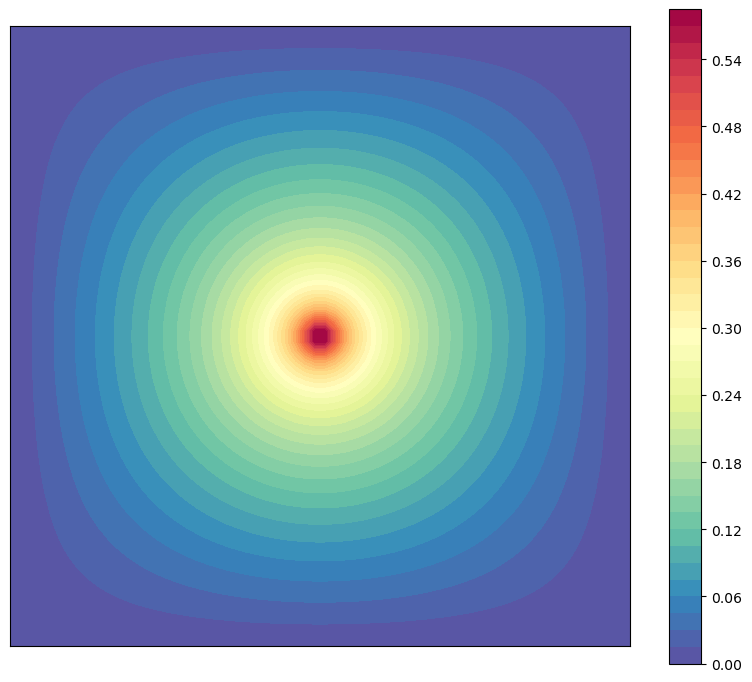} &
    \includegraphics[width=0.43\textwidth]{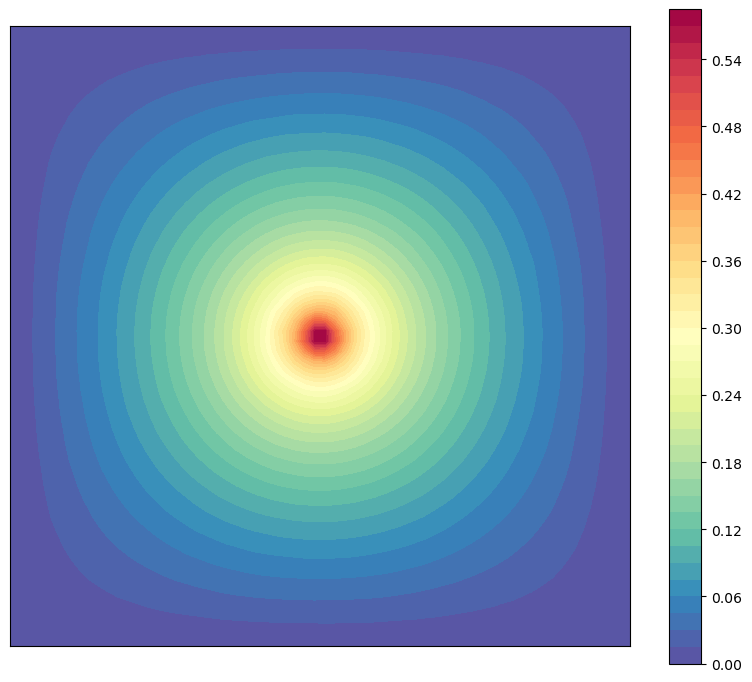} \\
    (a) Reference & (b) $\vT_\vxi^{(1)}$ \\
    \\
    $e_2 = 4.14 \times 10^{-3}$ & $e_2 = 2.23 \times 10^{-3}$\\
    \includegraphics[width=0.43\textwidth]{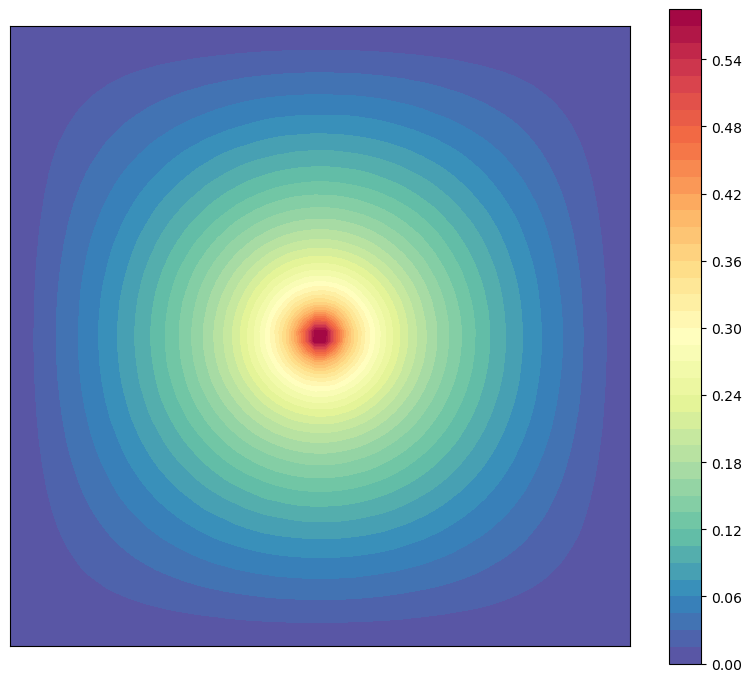} &
    \includegraphics[width=0.43\textwidth]{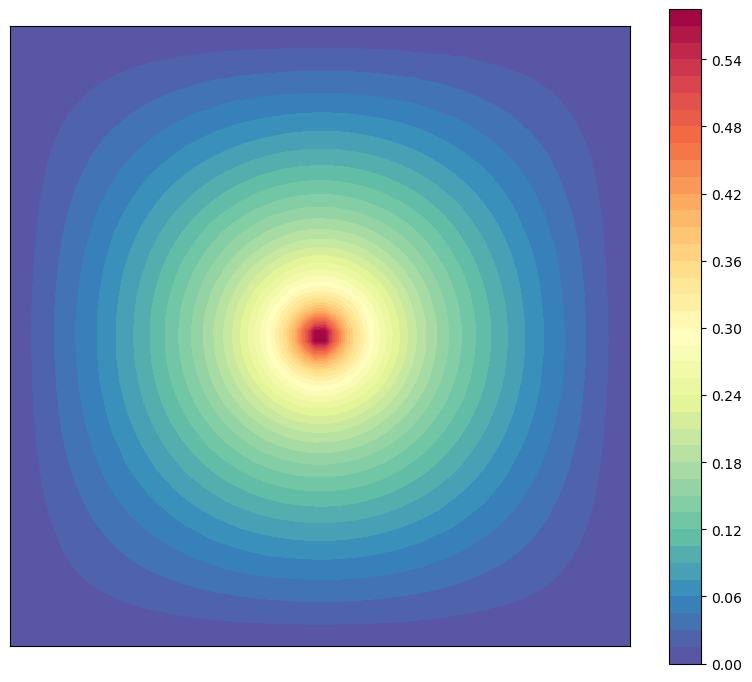} \\
    (c) $\vT_\vxi^{(2)}$ & (d) $\vT_\vxi^{(3)}$ 
  \end{tabular}
  \caption{Contour maps of the Green's functions $G(\vx,\vxi)$ at $\vxi=(0,0)$ for the reference solution and the predicted solutions  by the proposed deep surrogate model with different input forms. The corresponding $L^2$ errors (denoted as $e_2$) are also provided.}
  \label{fig:5}
\end{figure}

Here we investigate the effect of different input forms on the the learned Green's function. As mentioned in subsection \ref{sec:data_gen}, we provide three types of input tensors, including $1$-channel input $\vT_\vxi^{(1)}=\vrho_\vxi$,  $2$-channels input $\vT_\vxi^{(2)}=[\vR_\vxi, \vrho_\vxi]$, and  $3$-channels input $\vT_\vxi^{(3)} = [\vX, \vrho_\vxi]$.  
Figure \ref{fig:5} shows the contour maps of the Green's functions $G(\vx,\vxi)$ at $\vxi=(0,0)$ for the reference solution and the predicted solutions  by the proposed deep surrogate model with different input forms, where $Loss_{\text{jac}}$ is used.
We find that the results caused by these three input forms are almost the same, indicating that the performance of the deep surrogate model is mainly determined by the point source information and extra spatial location information isn't necessary. Therefore, we will use the first input form $\vT_{\vxi}^{(1)}$ in following experiments.

\subsubsection{Effect of the number of Jacobi  iterations in the constant strategy}\label{sec:Jacobi_fixed}
\begin{figure}[h!]
  \centering
  \begin{tabular}{ccc}
    $e_2 = 4.86\times 10^{-3}$& $e_2= 3.5\times 10^{-3}$\\
    \includegraphics[width=0.43\textwidth]{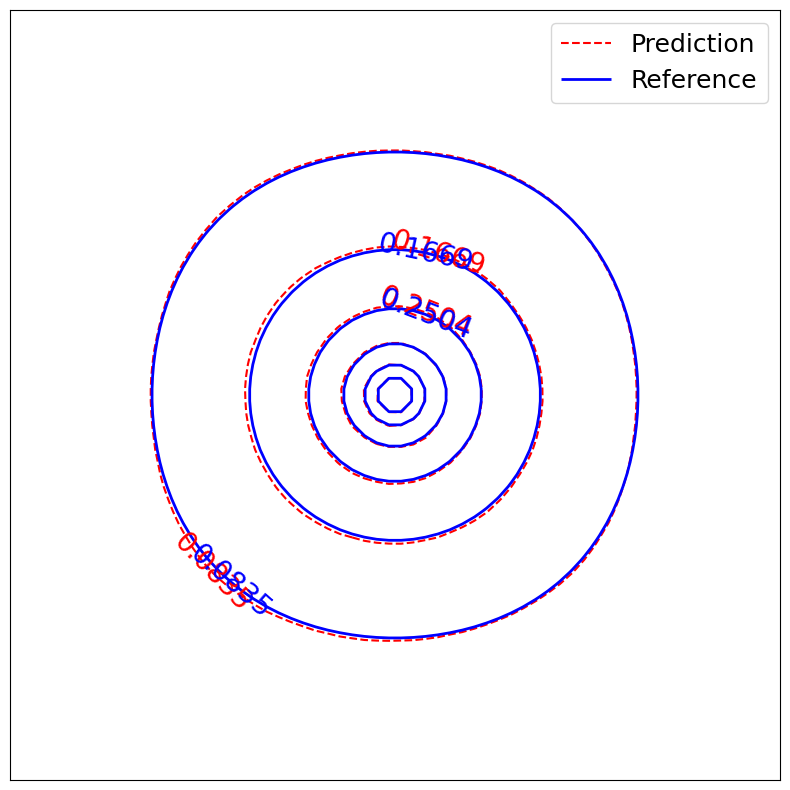} &
    \includegraphics[width=0.43\textwidth]{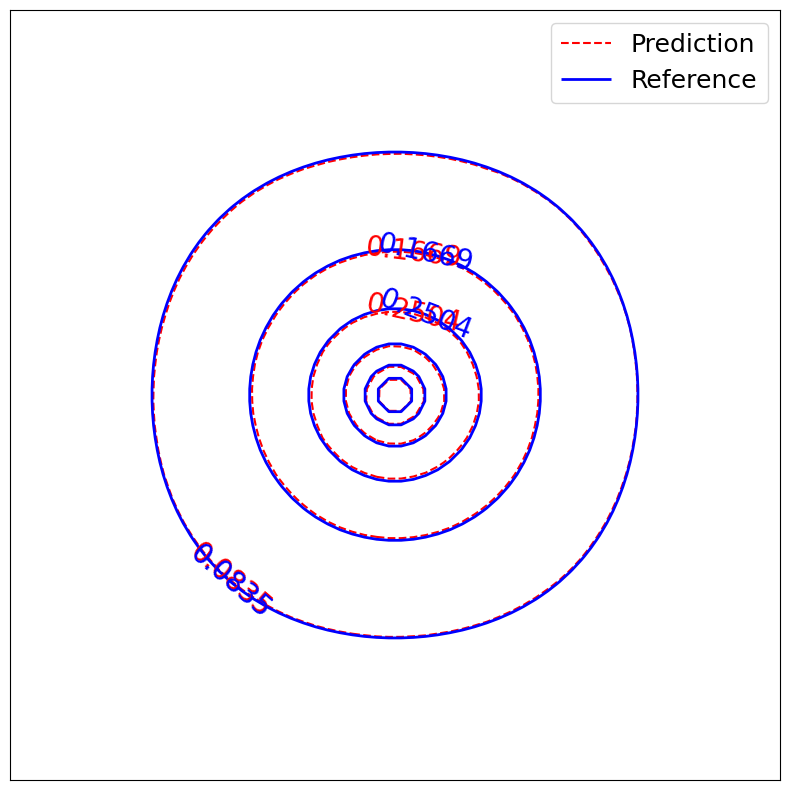} \\
    (a) $k=1$  & (b) $k=5$ \\
    \\
    \\
    $e_2 = 2.23\times 10^{-3}$& $e_2 = 4.72\times 10^{-3}$\\
    \includegraphics[width=0.43\textwidth]{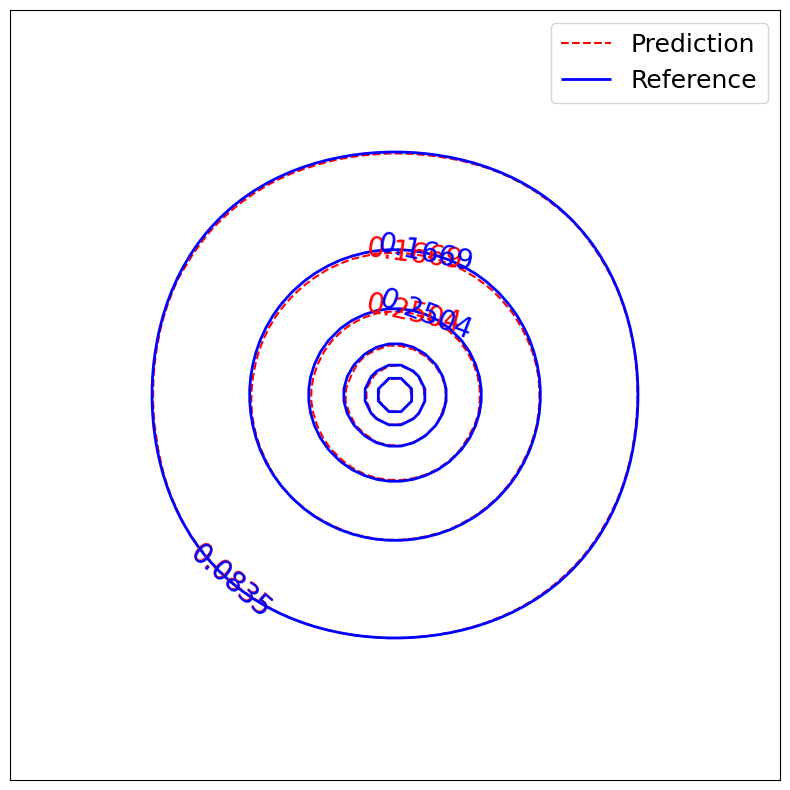} &
    \includegraphics[width=0.43\textwidth]{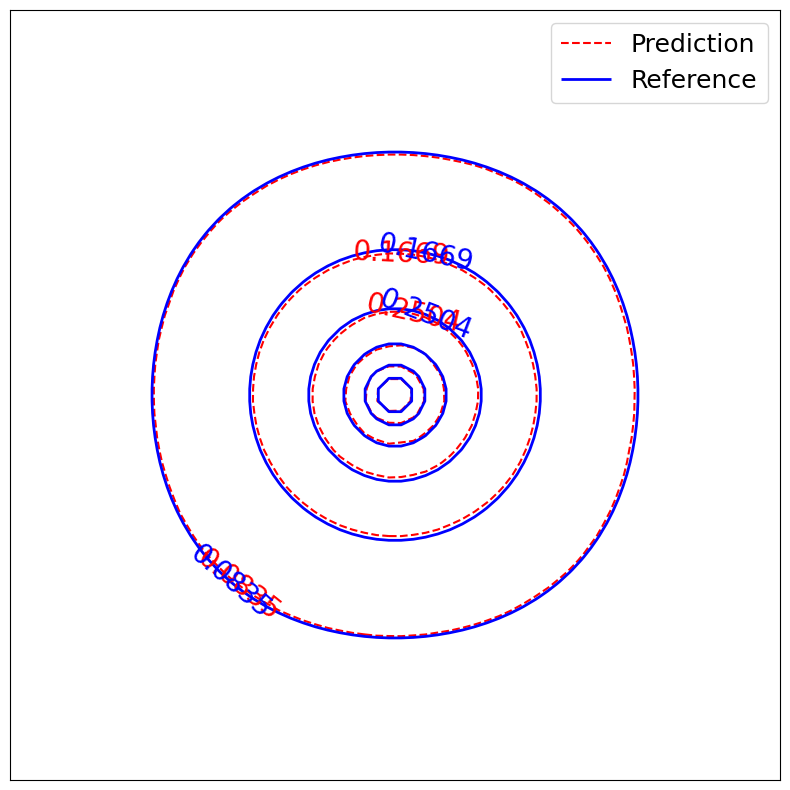} \\
    (c) $k=20$  & (d) $k=40$ 
  \end{tabular}
  \caption{Comparisons of contour maps of the Green's functions $G(\vx,\vxi)$ at $\vxi=(0,0)$ for the reference solution and the predicted solutions  by the proposed deep surrogate model with different number of Jacobi  iterations in  the constant strategy. The corresponding $L^2$ errors (denoted as $e_2$) are also provided.}
  \label{fig:Jacobi _iter}
\end{figure}

The Jacobi  scheme is a simple but important iterative method for solving large-scale linear systems. Here we study the effect of the number of Jacobi iterations in the constant strategy on the performance of the proposed deep surrogate model. Experimental results with a fixed number of Jacobi iterations $k$ are shown in Fig. \ref{fig:Jacobi _iter}. We find that $k$ has a significant impact on the model's performance: 1) when $k$ is set to be relatively small, the difference between the predicted solution and the reference solution is significant (see Fig. \ref{fig:Jacobi _iter}-(a) and  Fig. \ref{fig:Jacobi _iter}-(b)); 2)  as $k$ increases, the predicted solution gradually matches the reference solution (see Fig. \ref{fig:Jacobi _iter}-(c)), and subsequently, the number of iterations tends to be saturate, which means that further increase in $k$ may not  improve the predicted solution, and  may even lead to a poorer predictive performance of the model (see Fig. \ref{fig:Jacobi _iter}-(d)). In fact, if $k$ is set to be large enough, the approximate solution produced by the Jacobi iteration scheme is almost the exact solution, and our model  then could be regarded as the data-driven model. As mentioned earlier, the predicted solution generated by the data-driven surrogate model lacks some regularized constraints, which partially explains the phenomenon in Fig \ref{fig:Jacobi _iter}-(d).  That is to say, there is no need to choose large number for $k$ in practice. Although the incomplete Jacobi iterations may produce imperfect approximate solutions, it still provides a good estimate (label) for the training of the surrogate model, then this estimate is somehow corrected by the back-propagation algorithm. As the training progresses, a more accurate regularized solution will be generated in the end.

\subsubsection{Effect of he number of Jacobi  iterations in the dynamic and adaptive strategies}\label{sec:Jacobi_dynamic}



Now we  investigate the effect of the  other two training strategies, including the dynamic and adaptive ways to adjust the number of Jacobi iterations during the training process, on the model's performance. In our setting, the dynamic strategy is to initially set $k=40$, and then reduce $k$ by $10$ for every $20$ epochs until it is ultimately maintained at $10$. The adaptive strategy is to initially  set  $k=40$ for the first epoch,  and then let $k$ be adaptively adjusted in [0,20] the remaining process. Specifically, when $Loss_{\text{cur}} >1.2 Loss_{\text{pre}}$, multiply $k$ by $2$, and when $Loss_{\text{cur}} < 0.8 Loss_{\text{pre}}$, divide $k$ by $2$, where $Loss_{\text{cur}}$ and $Loss_{\text{pre}}$ are the validation losses at the current and previous epochs, respectively.

\begin{figure}[htbp]
  \centering
    \includegraphics[width=0.9\textwidth]{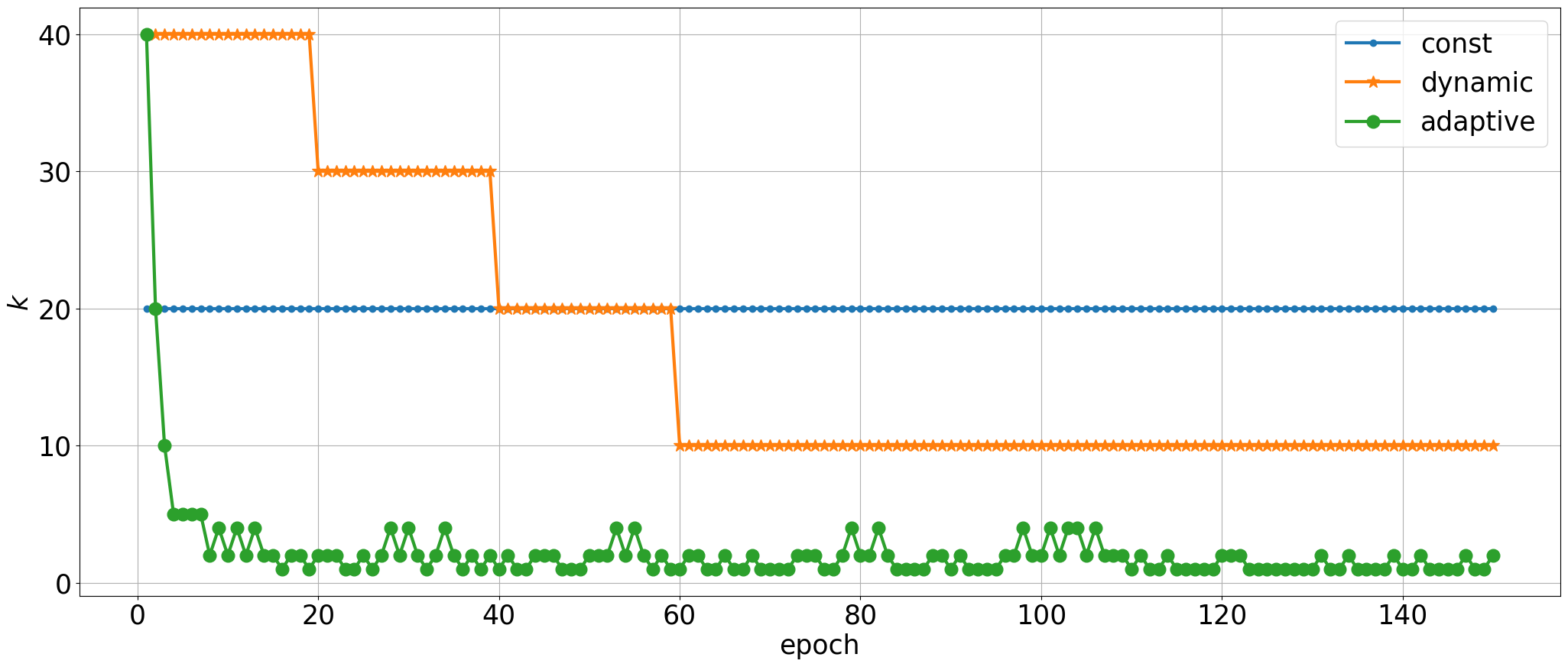} 
  \caption{Evolution of the number of Jacobi iterations  $k$ during the training process for three training strategies, where we fix $k=20$ for the constant strategy.}
  \label{fig:k}
\end{figure}

\begin{figure}[htbp]
  \centering
    \includegraphics[width=0.9\textwidth]{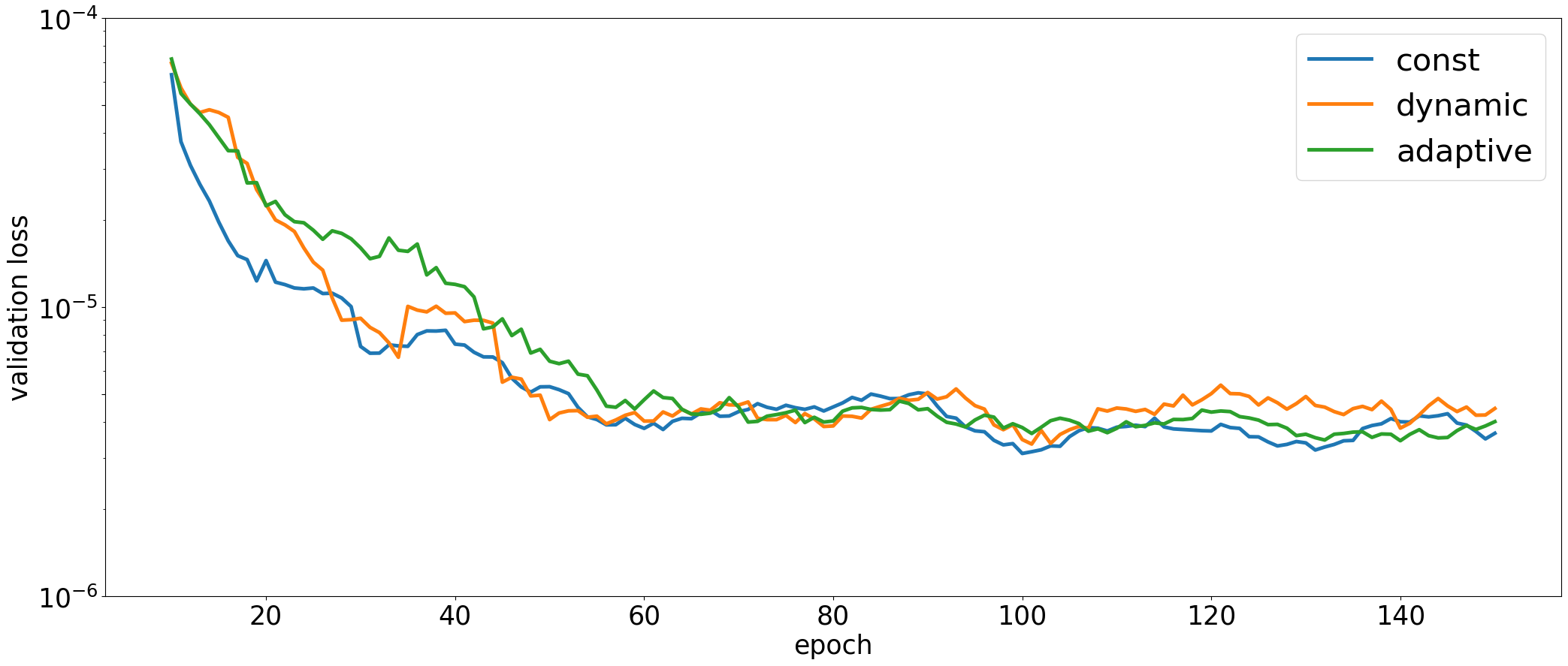} 
  \caption{Evolution of the validation loss during the training process for three training strategies. }
  \label{fig:loss}
\end{figure}

\begin{figure}[htbp]
  \centering
  \begin{tabular}{ccc}
    $e_2 = 2.23\times 10^{-3}$& $e_2 = 2.97\times 10^{-3}$ & $e_2 = 3.54 \times 10^{-3}$\\
   \hspace{-0.4cm} \includegraphics[width=0.32\textwidth]{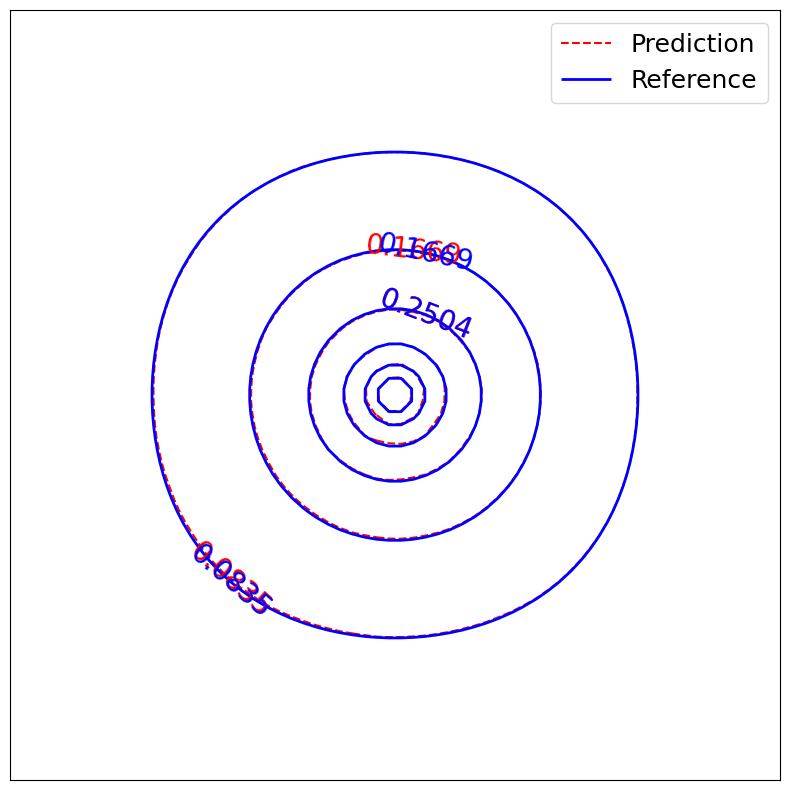} &
    \includegraphics[width=0.32\textwidth]{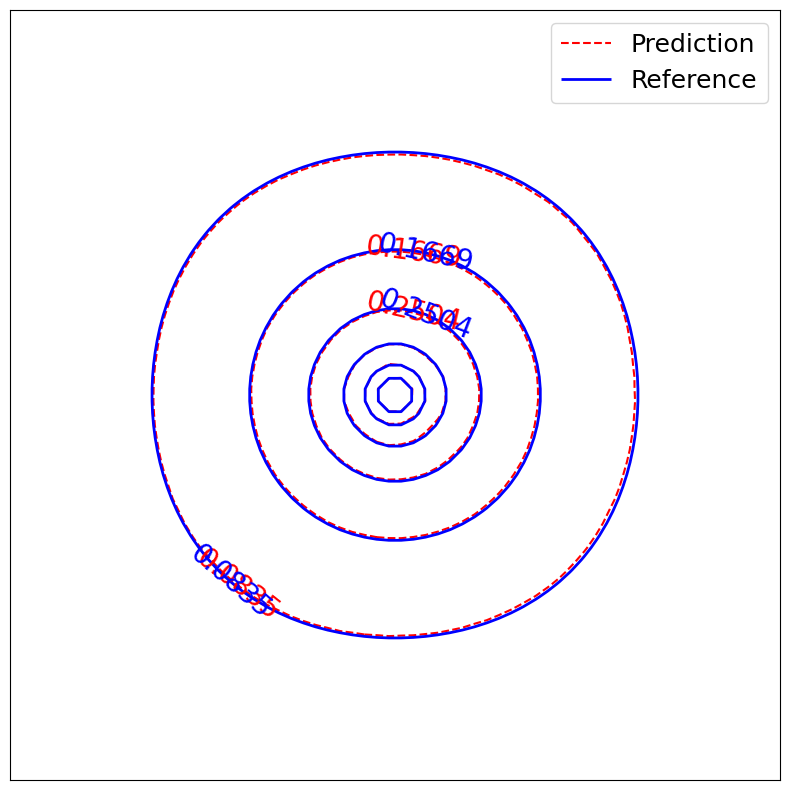} &
    \includegraphics[width=0.32\textwidth]{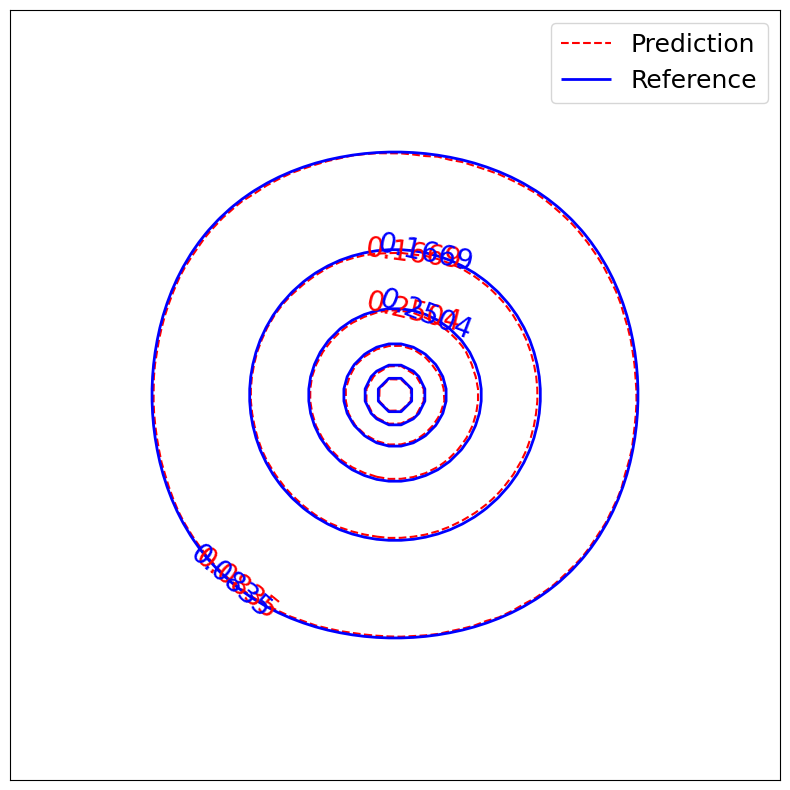} \\
    (a) constant  & (b) dynamic & (c) adaptive \\
  \end{tabular}
  \caption{
  Comparisons of contour maps of the Green's functions $G(\vx,\vxi)$ at $\vxi=(0,0)$ for the reference solution and the predicted solutions  by the proposed deep surrogate model with the three different training strategy. The corresponding $L^2$ errors (denoted as $e_2$) are also provided.}
  \label{fig:contour_Jacobi}
\end{figure}

Evolution of the number of Jacobi iterations  $k$ and the validation loss during the training process for three training strategies are shown in Fig. \ref{fig:k} and Fig. \ref{fig:loss}, respectively, where we fix $k=20$ for the constant strategy. Compared with the constant strategy, the other two strategies, especially the adaptive strategy, significantly reduce the total number of Jacobi iterations, while the downtrend of the loss function remains similar for all three strategies. We present the contour maps of the predicted Green's functions at $\vxi=(0,0)$ by the proposed deep surrogate model with the three different training strategies in Fig. \ref{fig:contour_Jacobi}, from which it can be seen that the dynamic and adaptive strategies significantly reduce computational complexity, while the predicted results are comparable to 
those of the constant strategy. Thus we will always use the adaptive strategy for all the remaining experiments.

\subsection{More examples for learning Green's functions}\label{sec:more_examples}
\begin{figure}[htbp]
  \centering
  \begin{tabular}{ccc}
  $e_2 = 2.23\times 10^{-3}$& $e_2 = 1.38\times 10^{-3}$& $e_2 =2.19\times 10^{-3}$\\
  \hspace{-0.4cm}   \includegraphics[width=0.32\textwidth]{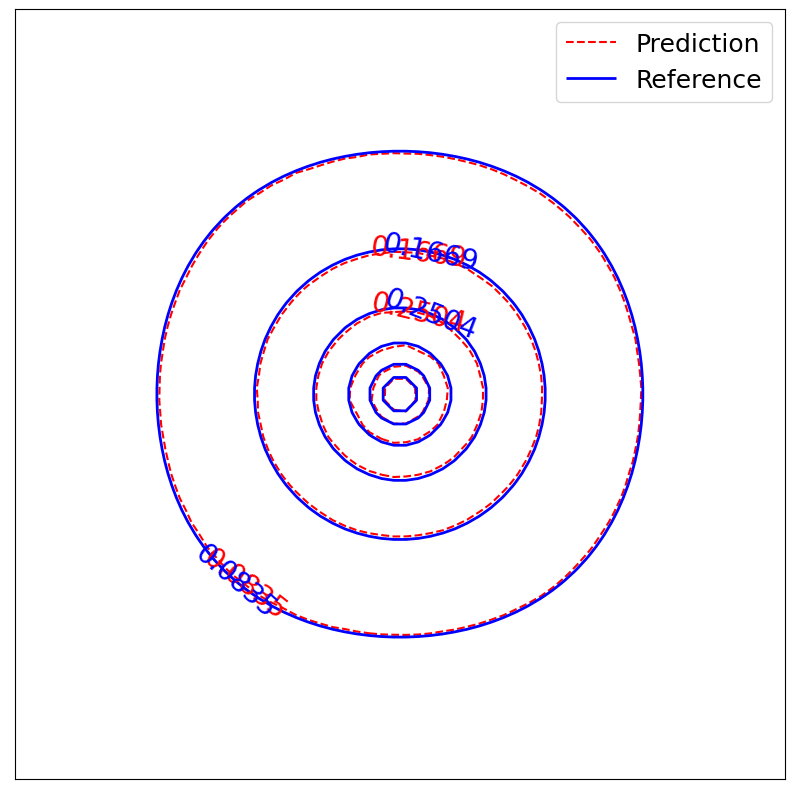} &
    \includegraphics[width=0.32\textwidth]{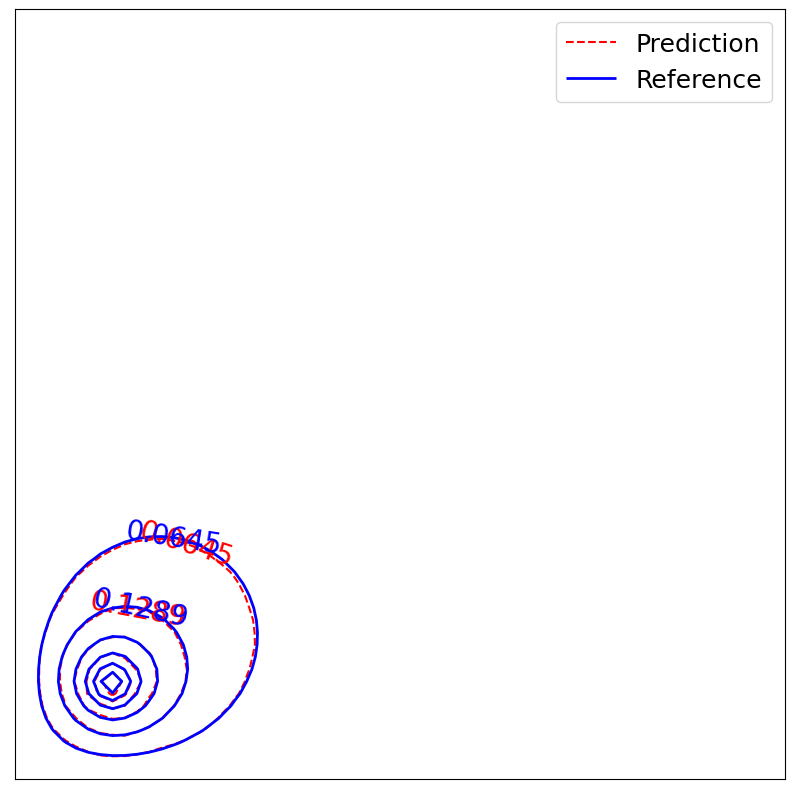} &
    \includegraphics[width=0.32\textwidth]{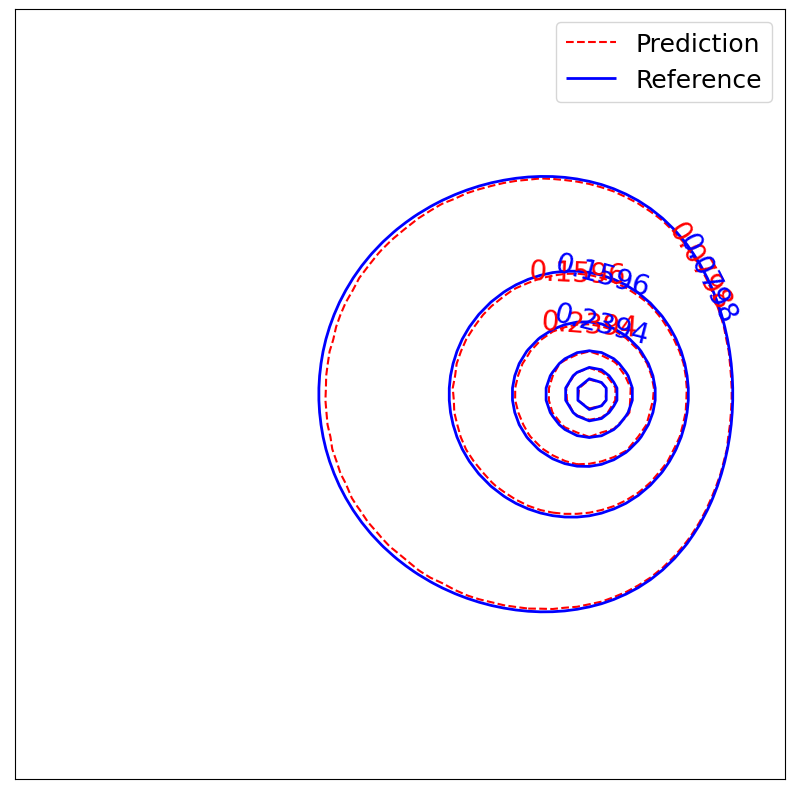} \\
    (a) $\vxi=(0,0)$ & (b) $\vxi=(-0.75,-0.75)$ & (c)  $\vxi=(0.5,0)$ 
  \end{tabular}
  \caption{Comparisons of contour maps of the Green's functions $G(\vx,\vxi)$ for the reference solution and the predicted solutions by the proposed deep surrogate model
   for the Laplacian operator with different source positions. The corresponding $L^2$ errors (denoted as $e_2$) are also provided.}
  \label{fig:Green_for_Laplacian}
\end{figure}

\begin{figure}[htbp]
  \centering
  \begin{tabular}{ccc}
    $e_2 = 3.11\times 10^{-3}$& $e_2 = 1.09\times 10^{-3}$& $e_2 = 2.44\times 10^{-3}$\\
\hspace{-0.4cm}     \includegraphics[width=0.32\textwidth]{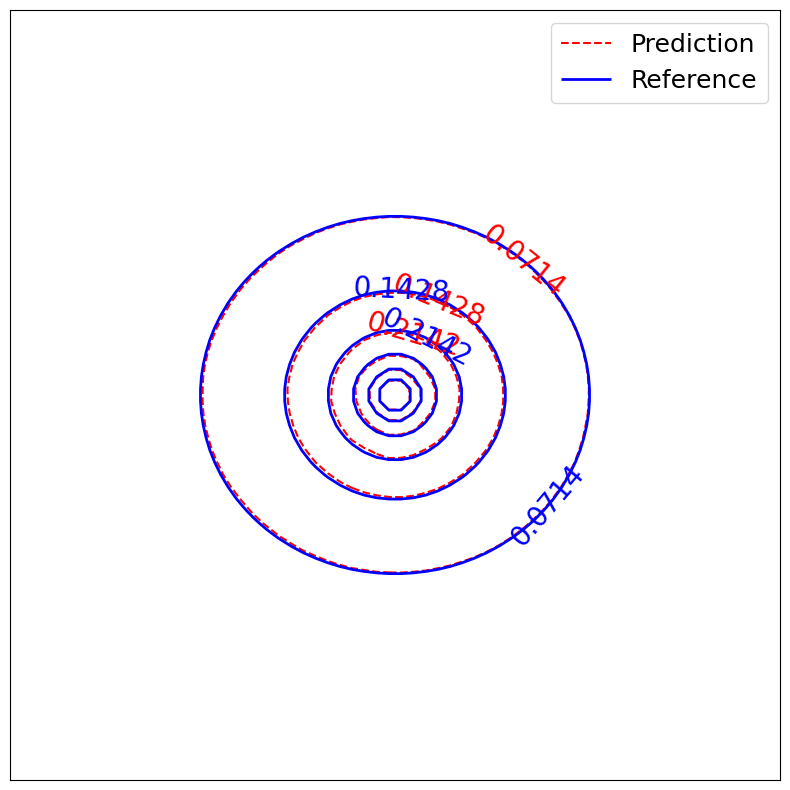} &
    \includegraphics[width=0.32\textwidth]{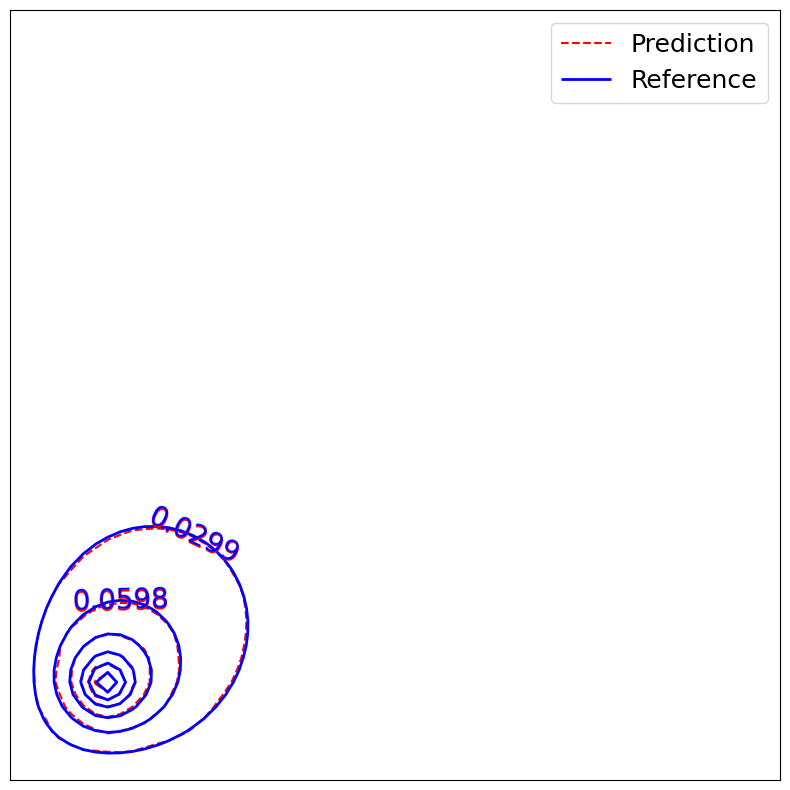} &
    \includegraphics[width=0.32\textwidth]{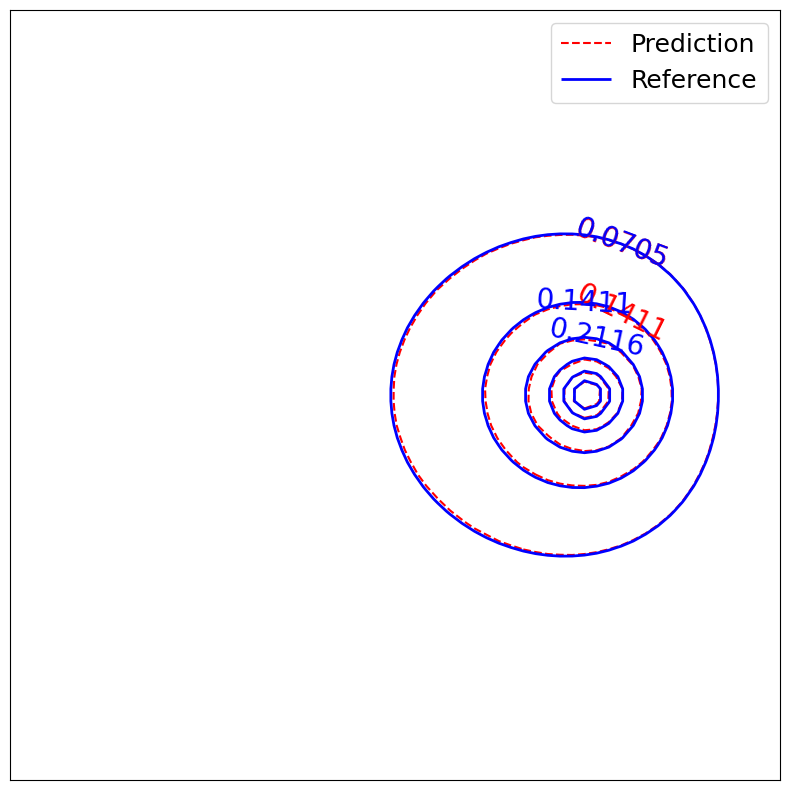} \\
    (a) $\vxi=(0,0)$  & (b) $\vxi=(-0.75,-0.75)$ & (c) $\vxi=(0.5,0)$
  \end{tabular}
  \caption{Comparisons of contour maps of the Green's functions $G(\vx,\vxi)$ for the reference solution and the predicted solutions by the proposed deep surrogate model
   for the reaction-diffusion operator \eqref{111} with different source positions. The corresponding $L^2$ errors (denoted as $e_2$) are also provided.}
  \label{fig:Green_for_RD}
\end{figure}

More experimental results about the learned Green's functions for different linear reaction-differential operators, produced by the proposed deep surrogate model, are provided in this subsection. Comparisons of contour maps (and corresponding numeral errors) shown in Fig. \ref{fig:Green_for_Laplacian} verify that our model can produce accurate prediction results for Green's function of Laplacian operator at different source positions, even near the boundary or corner of the domain.  We also investigate the deep surrogate model for learning Green's function of  the reaction-diffusion operator \eqref{eq:OP_RD} with the variable coefficients 
\begin{equation}\label{111}
a(\vx)=1+2x_2^2,\quad r(\vx)=1+x_1^2,
\end{equation}
 and the corresponding results are shown in Figure \ref{fig:Green_for_RD}, which demonstrate our  model again works very well..


\subsection{Fast solvers for solving PDEs}
In this subsection, we conduct some experiments for investigating performance of the fast solvers \eqref{eq:Quad}  based on the learned Green’s function in solving the model linear reaction-diffusion equation \eqref{eq:PDE_RD}. For the choice of numerical quadratures, we adopt $I_{\vx,h}^{R_l}$ as the 2D Simpson's rule  and $I_{\vx,h}^{E_m}$ as the 1D Simpson's rule, and the gradient $\nabla_\vx G(\vx,\vxi)$ on the boundary is approximated by a first order difference scheme. For the visualization and the computation of $L^2$ error between the exact solution and the approximate solution provided by fast solvers, we approximately compute $u(\vxi)$ by \eqref{eq:Quad} on a $64\times 64$ uniform mesh. 

\subsubsection{Laplacian equation}
\begin{figure}[htbp]
  \centering
  \begin{tabular}{ccc}
    & & $e_2=2.97\times 10^{-2}$ \\
  \hspace{-0.4cm}   \includegraphics[width=0.33\textwidth]{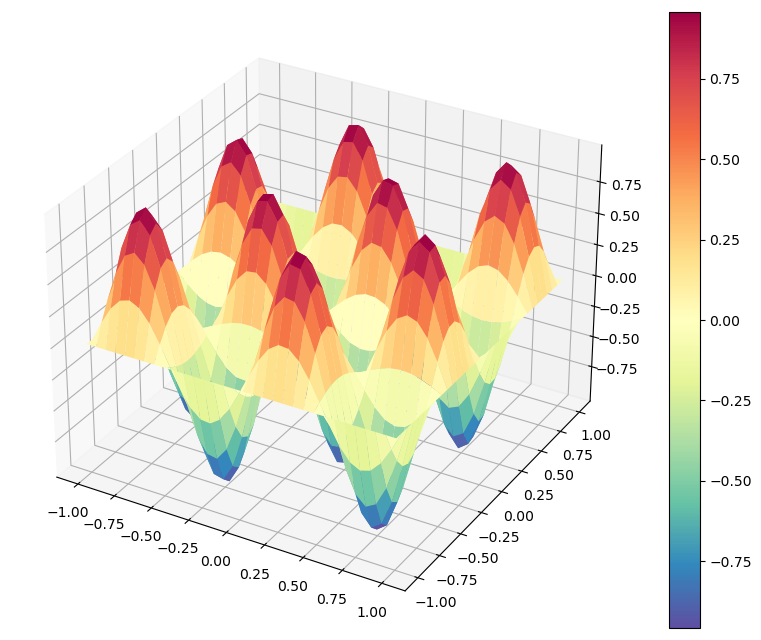} &
    \includegraphics[width=0.33\textwidth]{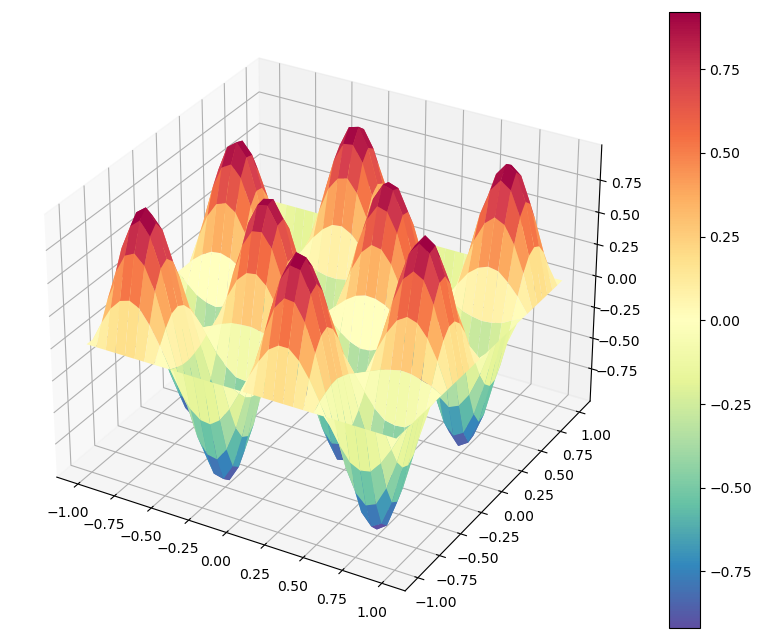} &
    \includegraphics[width=0.28\textwidth]{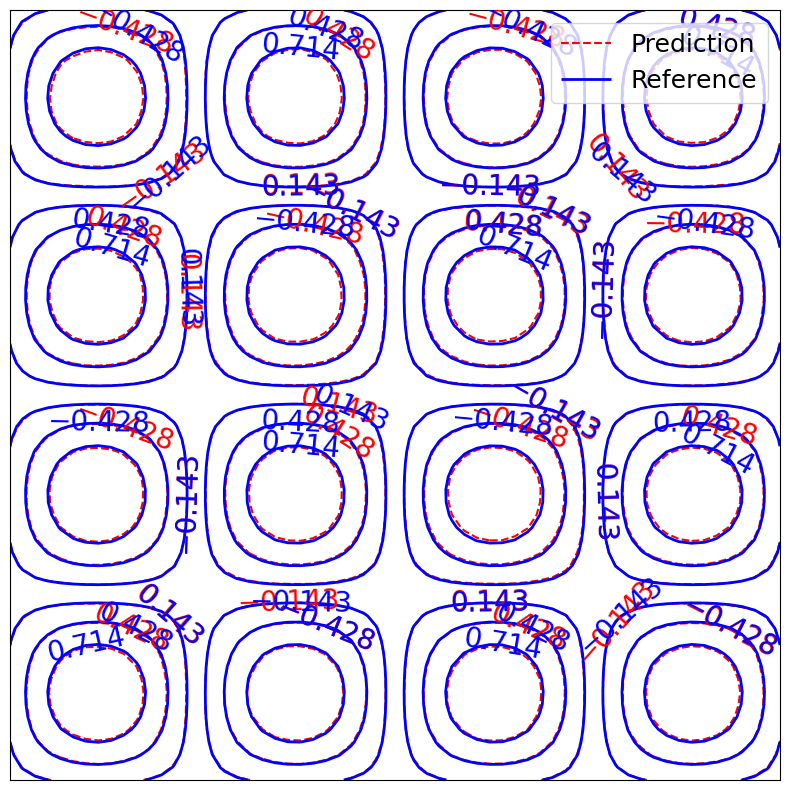} \\
    (a) reference & (b) prediction & (c) contour \\
  \end{tabular}
  \caption{Plots of the exact solution (left) and  the numerical solution (middle) produced by the fast solver based on the learned Green's function for the Poisson equation
  with  the solution \eqref{eq:lap_homo}  ($\lambda=2$). The comparison of their contour maps (right)  is also provided with the corresponding $L^2$ errors (denoted as $e_2$).}
  \label{fig:8}
\end{figure}

\begin{figure}[htbp]
  \centering
  \begin{tabular}{ccc}
     & & $e_2=5.98\times 10^{-2}$ \\
\hspace{-0.4cm}     \includegraphics[width=0.33\textwidth]{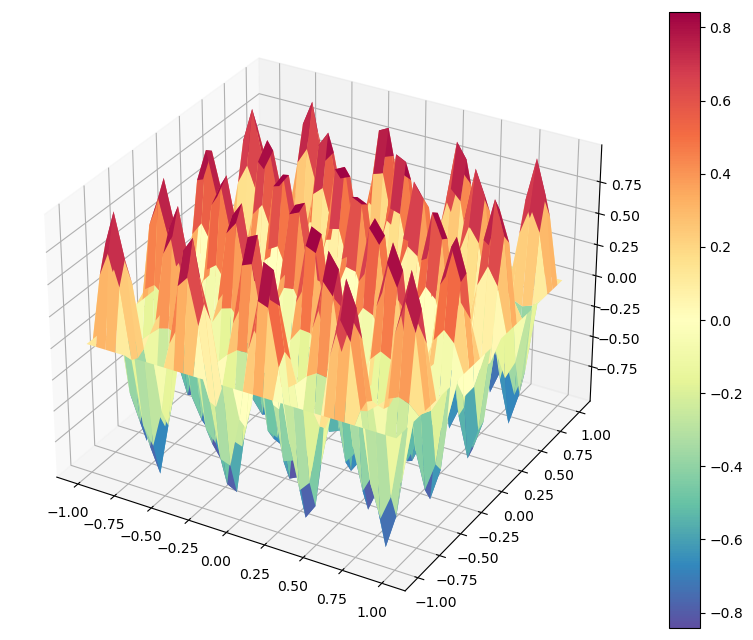} &
    \includegraphics[width=0.33\textwidth]{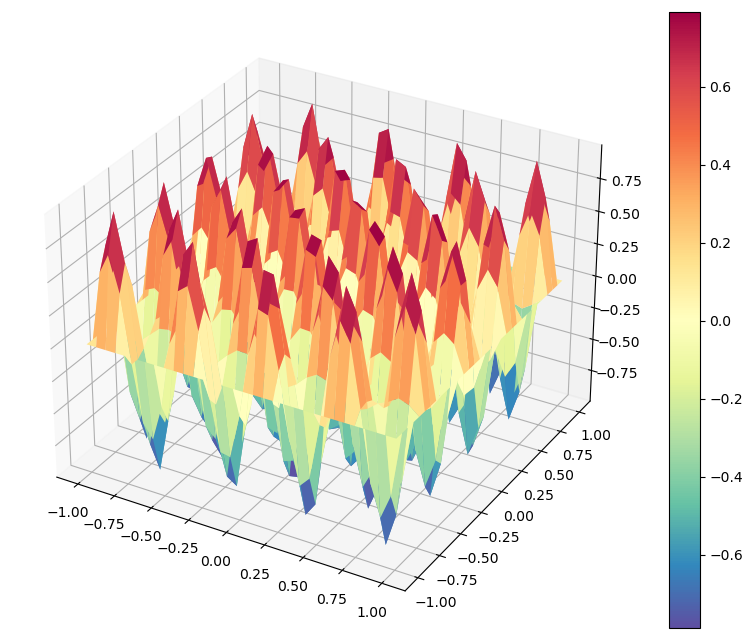} &
    \includegraphics[width=0.28\textwidth]{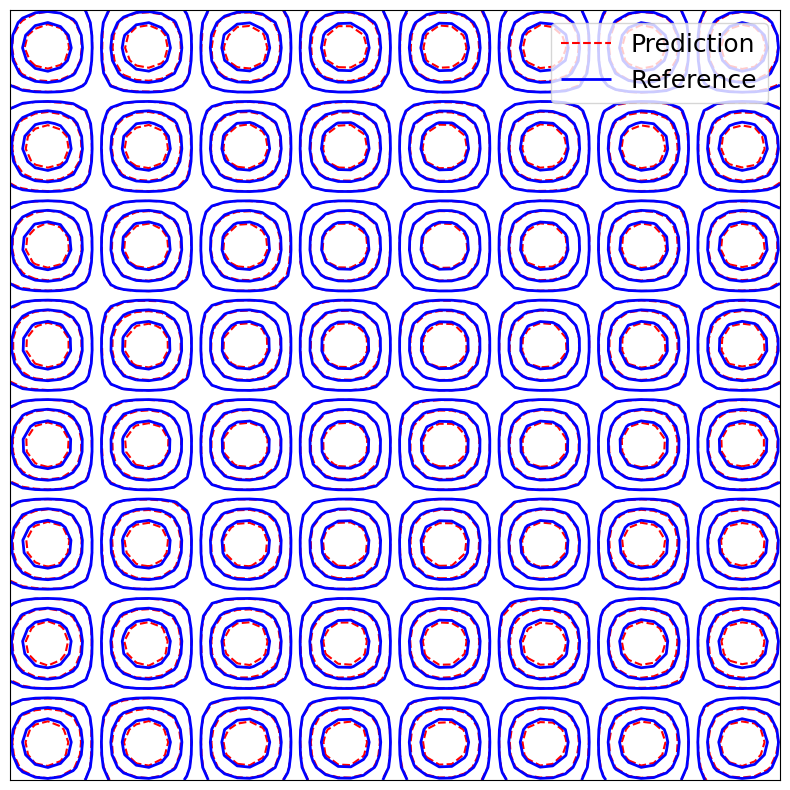} \\
    (a) reference & (b) prediction & (c) contour \\
  \end{tabular}
  \caption{
  Plots of the exact solution (left) and  the numerical solution (middle) produced by the fast solver based on the learned Green's function for the Poisson equation
  with  the solution \eqref{eq:lap_homo}  ($\lambda=4$). The comparison of their contour maps (right)  is also provided with the corresponding $L^2$ errors (denoted as $e_2$).}
  \label{fig:9}
\end{figure}

Let us consider the Laplacian equation.
First, we choose the exact solution as 
\begin{equation}\label{eq:lap_homo}
u(x_1,x_2)=\sin(2\lambda\pi x_1)\sin(2\lambda\pi x_2)
\end{equation}
 and the source term is determined accordingly, where $\lambda$ is used to determine the frequency of the solution. In this case, the boundary condition is homogenous, i.e., $g=0$, and the second integral in \eqref{eq:Green_representation} vanishes. 
 Plots of the exact solution, the numerical solution and their contour maps are shown  in Figure \ref{fig:8} for the case of  $\lambda=2$ and Figure \ref{fig:9} for the case of $\lambda=4$. We see that our fast solver produces good numerical solutions, and the $L^2$ errors  are within acceptable range.  

\begin{figure}[htbp]
  \centering
  \begin{tabular}{ccc}
     & & $e_2=3.12\times 10^{-2}$ \\
\hspace{-0.4cm}     \includegraphics[width=0.33\textwidth]{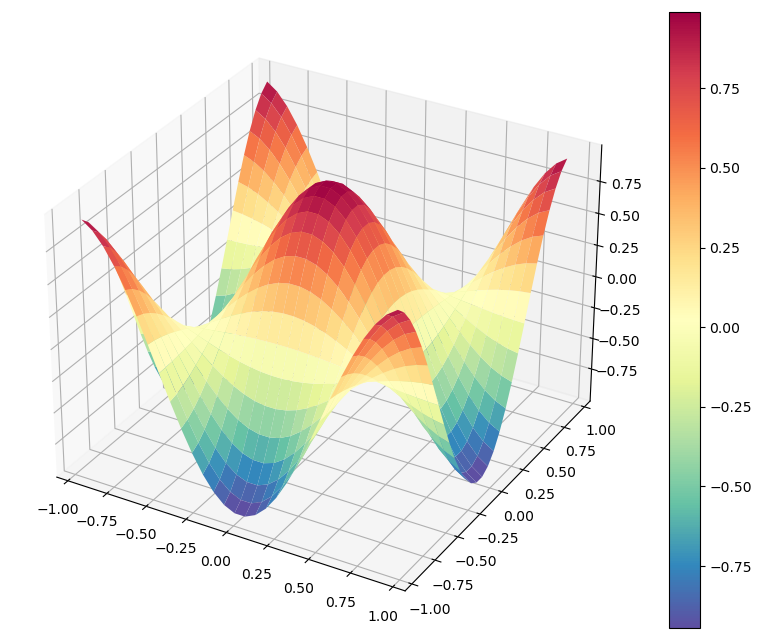} &
    \includegraphics[width=0.33\textwidth]{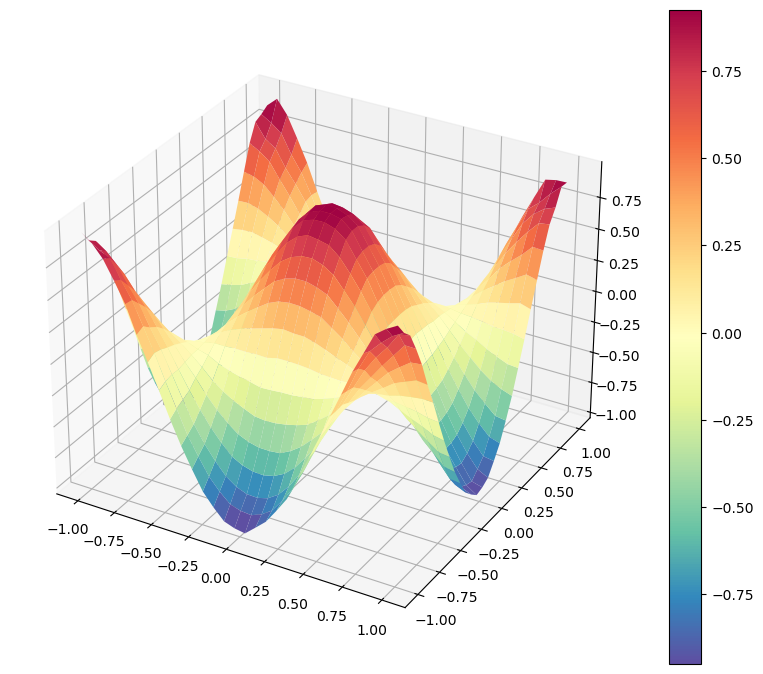} &
    \includegraphics[width=0.28\textwidth]{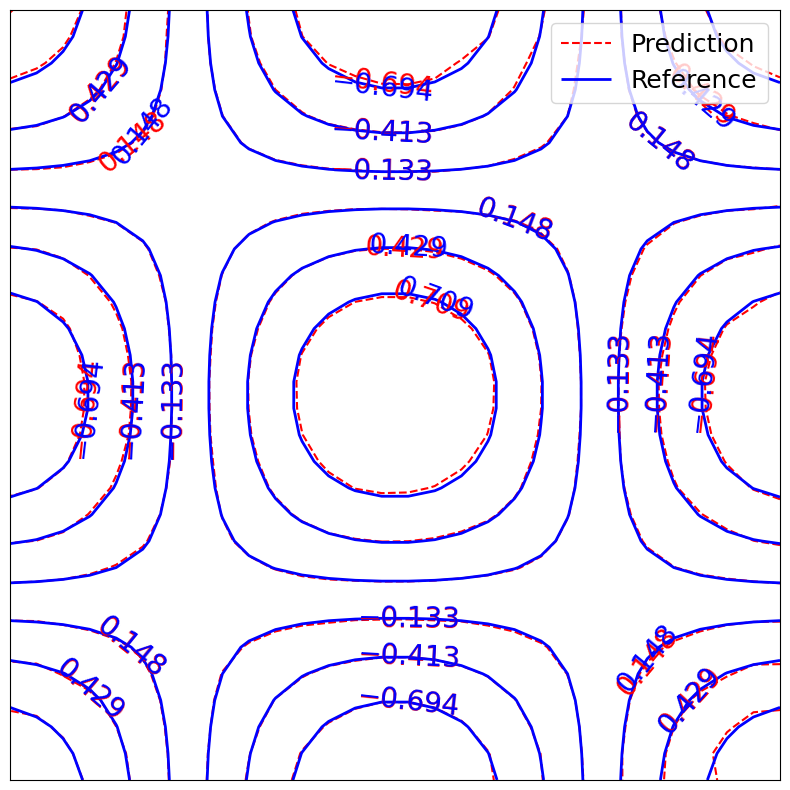} \\
    (a) reference & (b) prediction & (c) contour \\
  \end{tabular}
  \caption{
  Plots of the exact solution (left) and  the numerical solution (middle) produced by the fast solver based on the learned Green's function for the Poisson equation
  with  the solution \eqref{eq:lap_inhomo}. The comparison of their contour maps (right)  is also provided with the corresponding $L^2$ errors (denoted as $e_2$).}
  \label{fig:10}
\end{figure}

Second, we choose the exact solution as 
\begin{equation}\label{eq:lap_inhomo}
u(x_1,x_2)=\cos(\pi x_1)\cos(\pi x_2)
\end{equation}
 and the boundary condition (inhomogeneous now) and source term are determined accordingly.  Plots of the exact solution, the numerical solution and their contour maps are shown in  Figure \ref{fig:10}, from which we see that the fast solver  is also suitable for solving the Laplacian equation with inhomogeneous boundary condition. 
 

\subsubsection{The reaction-diffusion equation}
\begin{figure}[htbp]
  \centering
  \begin{tabular}{ccc}
    & & $e_2=8.31\times 10^{-3}$ \\
\hspace{-0.4cm}     \includegraphics[width=0.33\textwidth]{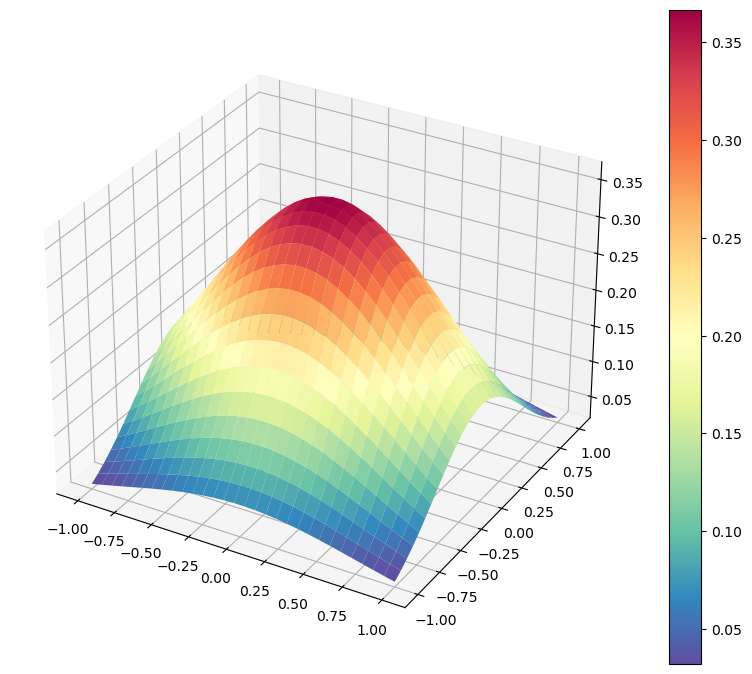} &
    \includegraphics[width=0.33\textwidth]{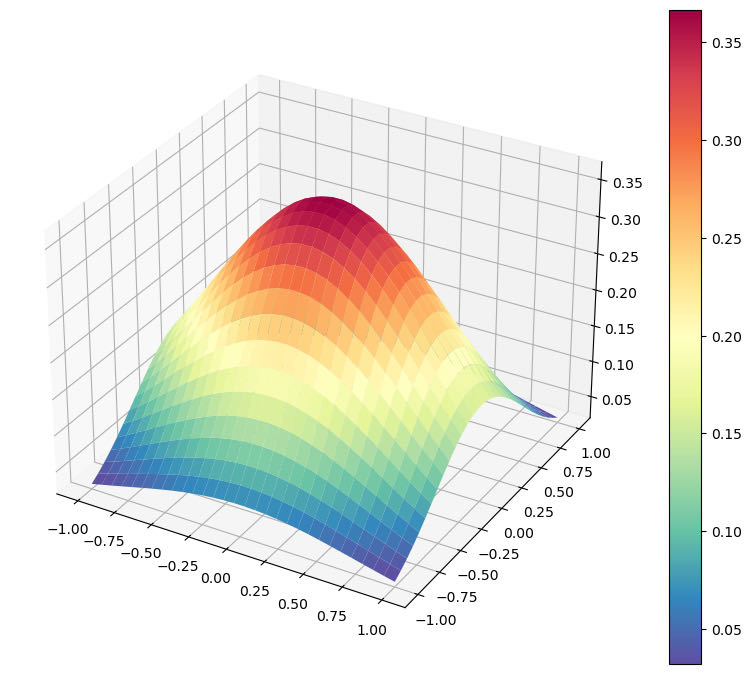} &
    \includegraphics[width=0.28\textwidth]{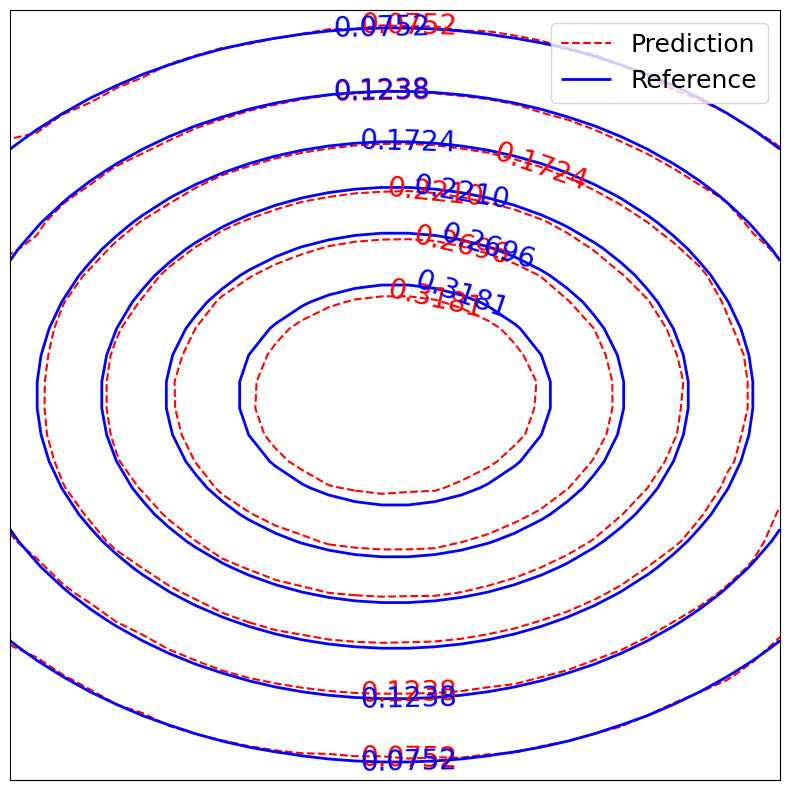} \\
    (a) reference & (b) prediction & (c) contour \\

  \end{tabular}
  \caption{Plots of the exact solution (left) and  the numerical solution (middle) produced by the fast solver based on the learned Green's function for the linear reaction-diffusion equation
  with the coefficients \eqref{111} and the solution \eqref{eq:reaction_inhomo}. The comparison of their contour maps (right)  is also provided with  the corresponding $L^2$ errors (denoted as $e_2$).}
  \label{fig:11}
\end{figure}

We consider the reaction-diffusion equation with the coefficients defined in \eqref{111} and the exact solution is set to be
\begin{equation}\label{eq:reaction_inhomo}
u(x_1,x_2)=10^{-(x_1^2+2x_2^2+1)}
\end{equation}
 The boundary condition and source term are then determined accordingly. Plots of the exact solution, the numerical solution and their contour maps are shown  in Figure \ref{fig:11}, from which we see that numerical solution produced by the fast solver \eqref{eq:Quad} based on the learned Green's function again matches the exact solution very well. 

\section{Conclusion}

In this paper we propose and numerically study a deep surrogate model for learning Green's function of linear reaction-diffusion operator based on the U-Net architecture. Inspired by the Jacobi  iteration scheme for solving linear systems, a novel Jacobi-type loss function and corresponding training strategies are designed and demonstrated to be very effective. 
In addition, a fast solver is tested and shown to be effective  for numerical solution of linear reaction-diffusion equations based on the learned Green's function. 
The proposed model is a beneficial attempt to integrate deep learning with traditional numerical methods. It fully utilizes the powerful expression capability of neural networks, and on the other hand, it also combines advantages of traditional numerical methods.



\bibliographystyle{unsrt}
\bibliography{main}







\end{document}